\documentstyle[12pt]{article}
\def\BP{\vskip 10pt}

\def\bss{\backslash}
\def\F{\noindent}
\def\newpage{\vfill\eject}
\def\q{\quad\quad}
\def\qq{\quad\quad\quad\quad}

\def\SP{\vskip 5pt}

\def\BP{\vskip 10pt}

\def\split{\begin{array}{ll}}
\def\endsplit{\end{array}}
\def\cen{\centerline}
\def\dm{\displaystyle}
\def\mb{\mbox}
\def\r{\,$\mb}

\def\A{\;\,$}

\def\bB{{\bf B}}
\def\bC{{\bf C}}

\def\bR{{\bf R}}

\def\CD{{\cal D}}

\def\Con{C_0^\infty}

\def\f{$\,}
\def\gr{\nabla}

\def\m{\vert}
\def\M{\Vert}
\def\pa{\partial}

\def\Ro{{\bf R}}
\def\Rb{{\bf R}^2}
\def\R3{{\bf R}^3}

\def\RN{{\bf R}^N}

\def\Sch{Schr\"odinger}

\def\T{\eqno}

\def\tx{{\widetilde x}}

\def\Ag{1}
\def\BDG{2}
\def\BM{3}
\def\DPa{4}
\def\DPb{5}
\def\Ei{6}
\def\IS{7}
\def\JS{8}
\def\JSa{9}
\def\RZ{10}
\def\Sa{11}
\def\Sc{12}
\def\Sd{13}
\def\Wea{14}
\def\Web{15}
\def\Wi{16}
\def\Z{17}

\begin{document}

$$
\ \ 
$$

\vskip 60pt

\centerline{\bf  On the Spectrum of the Reduced Wave Operator}
\vskip 5pt
\centerline{\bf  with Cylindrical Discontinuity$^*$}
\vskip 30pt
\centerline{Willi J\"{a}ger}
\vskip 10pt
\centerline{Department of Mathematics} 
\centerline{The University of Heidelberg}
\centerline{D-69120 Heidelberg 1}
\centerline{Germany}
\vskip 10pt
\cen{and}
\vskip 10pt
\centerline{Yoshimi Sait$\bar { \hbox{o} }$}
\vskip 10pt
\centerline{Department of Mathematics} 
\centerline{University of Alabama at Birmingham}
\centerline{Birmingham, Alabama 35294}
\centerline{U. S. A.}
%
%            The end of Topmatter
\vskip 160pt

       *) This work was supported by Deutche Forschungs Gemeinschaft
      through SFB 359.

\newpage

\cen{  {\bf \S1. Introduction.} }              
\BP

       Consider the differential expression
$$
             h = -\,\mu(x)^{-1}\Delta.                       \T (1.1)
$$
     Here \f \Delta \A is the Laplacian in \f \RN \A with \f N \ge 2 $, and
     \f \mu(x) \A is a positive function on \f \RN \A given by
$$
        \mu(x) =
\left\{ \split
           \mu_1 \qq (x \in \Omega_1), \\                                   
           \mu_2 \qq (x \in \Omega_2), \\                                   
\endsplit \right.                                      \T (1.2)
$$
     where \f \mu_1, \mu_2 > 0$, \f \mu_1 \ne \mu_2$, and \f \Omega_{\ell} $, 
     \f \ell = 1, 2 $, are open sets of \f \RN \A such that
$$      
\left\{ \split
        \Omega_1 \cap \Omega_2 = \emptyset, \\
        \overline{\Omega_1} \cup \Omega_2 
                = \Omega_1 \cup \overline{\Omega_2} = \RN,
\endsplit \right.                                       \T (1.3)
$$
     \f \overline{\Omega_{\ell}} \A being the closure of
     \f \Omega_{\ell} $. It is easy to see that a selfadjoint realization 
     \f H \A of \f h \A is given by
$$
\left\{ \split
            D(H) = H^2(\RN), \\
            Hu = hu,
\endsplit \right.                                       \T (1.4)           
$$
     in the Hilbert space 
$$
               X = L_2(\RN; \mu(x)dx),                   \T (1.5)
$$
     where \f D(H) \A denotes the domain of \f H$, \f H^k(\RN) \A denotes 
     the \f k $-th order Sobolev space over \f \RN$, and \f hu \A should be 
     taken in the sense of distributions. The separating surface \f S \A 
     is defined by
$$
    S = \overline{\Omega_1} \cap \overline{\Omega_2} 
                               = \pa\Omega_1 = \pa\Omega_2,   \T (1.6)
$$
     \f \pa\Omega_{\ell} $, \f \ell = 1, 2 $, being the boundary of
     \f \Omega_{\ell} $. Let 
$$     
      n^{(\ell)}(x) = (n_1^{(\ell)}(x), n_2^{(\ell)}(x),     
                           \cdots, n_N^{(\ell)}(x)) \qq (\ell = 1, 2)
                                                              \T (1.7)
$$     
     be the unit outward normal of \f \Omega_{\ell} \A at \f x \in 
     \pa\Omega_{\ell} = S$. Obviously we have \f n^{(1)}(x) + n^{(2)}(x)
     = 0 \A for \f x \in S$. Eidus\,[\Ei] considered the operator \f H \A 
     under the following assumptions:\ {\it there exist positive constants 
     \f c_1 \A and \f c_2 \A such that
$$
            \big| n_N^{(1)}(x) \big| \ge c_1   \qq (x \in S),        \T (1.8)
$$
     and     
$$
        \m x \cdot n^{(1)}(x) \m \le c_2  \qq (x \in S),             \T (1.9)
$$
     where \f x \cdot n^{(1)}(x) \A is the inner product of \f x \A and
     \f n^{(1)}(x) \A in \f \RN $.}\ Note that a cone having its vertex at 
     the origin and the positive \f x_N$-axis as its axis satisfies (1.8) and 
     (1.9). Imposing the above assumptions, Eidus\,[\Ei] proved the limiting 
     absorption principle for \f H $, that is, by denoting by \f R(z) \A 
     the resolvent of \f H $, the limits
$$
          \lim_{\eta\downarrow 0} R(\lambda \pm i\eta) = R_{\pm}(\lambda)
                     \qq {\rm in\ \ } \bB(L_{2, 1}(\RN), L_{2, -1}(\RN))
                                                                 \T (1.10)       
$$
     exist for \f \lambda > 0$, where the weighted \f L_2 \A space 
     \f L_{2, t}(\RN)$, \f t \in \bR $, is defined by
$$
      L_{2, t}(\RN) = \{ f \ :\ (1 + \m x \m)^{t}f(x) \in L_2(\RN) \},
                                                                  \T (1.11)         
$$
     and \f \bB(X, Y) \A is the Banach space of all bounded linear operators
     from \f X \A into \f Y $. Then, Sait\={o}\,[\Sd] showed that 
     \f L_{2,1}(\RN) \A and \f L_{2,-1}(\RN) \A in (1.10) can be replaced by     
     \f L_{2, \delta}(\RN) \A and \f L_{2, -\delta}(\RN) \A with 
     \f \delta > 1/2 $, respectively. This means that the limiting 
     absorption principle for \f H \A holds on the same weighted \f L_2 \A 
     spaces as are used for the \Sch\ operator (cf. Agmon\,[\Ag], 
     Ikebe-Sait\={o}\,[\IS] and Sait\={o}\,[\Sa]). Then Roach-Zhang\,[\RZ]
     has shown that \f u = R^{\pm}(\lambda)f $, where \f \lambda > 0 \A and
     \f f \in L_{2, \delta}(\RN) \A with \f \delta > 1/2$, is characterized 
     as a unique solution of the equation
$$
             (-\mu(x)^{-1}\Delta - \lambda)u = f             \T (1.12)
$$
     with the radiation condition
$$
   \lim_{R\to\infty} 
       \frac1R \int_{B_R} |\gr u \mp i\sqrt{\lambda\mu(x)}
                         \tx u|^2\, dx = 0  \qq   (\tx = \frac{x}{|x|}),
                                                                \T (1.13)
$$
     \f B_R \A being the ball with radius \f R \A and center at the origin.
     The condition (1.13) is a natural extension of the radiation condition 
     for the \Sch\ operators ([\IS], [\Sa]). [\RZ] also gave another proof 
     of the limiting absorption principle for \f H $.

        In this work we are going to show the limiting absorption principle
     for \f H \A whose separating surface \f S \A satisfies a new
     condition (see Assumption 2.1) so that we can treat, for example, 
     the case where \f \Omega_1 \A is an infinite cylindrical domain. 
     Our proof of the limiting absorption principle will show that not 
     only the uniqueness of the solution but also the existence of the 
     limit (1.10) can be proved through the estimate of the radiation 
     condition term 
$$
            \CD u = \gr u + \{(N-1)/(2r)\}\tx u - ik\tx u,    \T (1.14)
$$
     where 
$$
\left\{ \split
         u = R(z)f,  \\
         f \in X, \\
         k = k(x, z) = \sqrt{z\mu(x)}, \\
         r = |x|,    \\
         \tx = x/|x|.
\endsplit \right.                                            \T (1.15)         
$$

        The method demonstrated here can be applied to some other cases where the
     number of the media is more than \f 2 \A or infinitely many (multimedia 
     cases). We shall discuss these cases with its short-range or long-range
     perturbation eleswhere ([\JSa]).

        Another multimedia problem which has been discussed extensively
     is the stratified media in which the coefficients of the operator
     are the functions of \f x' \in {\bf R}^k \subset \RN $, \f k < N$.
     Some pertubed operators of the above type have been discussed, too.
     Here we refer Wilcox\,[\Wi], Ben-Artzi-Dermanjian-Guillot\,[\BDG], 
     Weder\,[\Wea], [\Web], DeBi\'{e}vre-Pravica\,[\DPa], [\DPb], 
     Boutet de Monvel-Berthier-Manda [\BM], and Zhang\,[\Z]. 
     In [\DPb] S. DeBi\'{e}vre and D. W. Pravica proved that there is no 
     point spectrum for the stratified propagators without any additional 
     conditions other than sufficient smoothness of the coefficients 
     at infinity. This is an extention of R. Weder\, [\Wea]. In this
     work and also in the work [\Ei], [\Sd], and [\RZ], we are interested
     in the non stratified case, in which it seems that the absence of the 
     point spectrum can not be obtained without imposing some additional 
     conditions.

        In \S2 we introduce the conditions on the separating surface \f S \A 
     and the function \f \mu(x) $. In \S3 $\sim$ \S6 we assume that 
     \f N \ge 3 $. The uniqueness of the solution of the equation (1.12) 
     with
$$
         \liminf_{R\to\infty} \int_{S_R} |\CD u|^2 \, dS = 0,   
                                                                 \T (1.16)
$$
     where \f S_R \A is the sphere with radius \f R \A and center at the 
     origin, will be shown in \S3. Our starting point in \S3 is an identity
     involving the radiation condition term \f \CD u \A (Proposition 3.3).
     This is an extension of a similar identity in the 
     case of \Sch\ operator ([\IS], [\Sc]). Proposition 3.3 is also used 
     in \S4, where an estimate for the radiation condition term \f \CD u \A 
     is given. In \S5 some more estimates for \f u = R(z)f \A will be given, 
     and these estimates are combined in \S6 to give the proof of the 
     limiting absorption principle for \f H $. We discuss the case that 
     \f N = 2 \A in \S7 since we treat this case in a slightly different
     way although the result is rather similar to the case of \f N \ge 3 $.
     As for some technichal details of the computations and arguments 
     appeared In \S3 $\sim$ \S7, we refer to 
     J\"{a}ger-Sait$\bar{\hbox{o}}$\,[\JS].
\SP

        {\bf Acknowledgement.} This work was finished when the second author
     was visiting the University of Heidelberg from October 1994 through 
     March 1995. Here he would like to thank Deutsche Forschungs Gemeinschaft for
     its support through SFB 359. Also the second author is thankful to 
     Professor Willi J\"{a}ger for his kind hospitality during this period.

$$
\ \ \ \ \ 
$$

\cen{     {\bf \S2. The operator \f H = -\, \mu(x)^{-1}\Delta $.}}   
\BP

        We shall start with describing the conditions imposed on our 
     operator.
\BP

        {\bf Assumption 2.1.}\ (i) Let \f N \A be a positive integer such 
     that \f N \ge 2 \A and let \f \Omega_{\ell} $, \f \ell = 1, 2 $, are 
     open sets of \f \RN \A satisfying (1.3). 
\SP
     
        (ii) Let the separating surface \f S \A be defined by (1.6). The 
     separating surface \f S \A is assumed to be 
     an $N-1$-dimensional continuous surface which consists of a finite
     number of smooth surfaces.     
\SP
     
        (iii) Let \f \mu(x) \A is a positive function on \f \RN \A given by
$$
        \mu(x) =
\left\{ \split
           \mu_1 \qq (x \in \Omega_1), \\  *[4pt]                                 
           \mu_2 \qq (x \in \Omega_2), \\                                   
\endsplit \right.                                       \T (2.1)
$$
     where \f \mu_1, \mu_2 > 0$, \f \mu_1 \ne \mu_2$. Further, we assume 
     that
$$
       (\mu_2 - \mu_1)(x \cdot n^{(1)})
          = (\mu_1 - \mu_2)(x \cdot n^{(2)}) \ge 0        \T (2.2)
$$
     for almost all \f x \in S $, where \f n^{(\ell)} $, \f \ell = 1, 2 $,
     is the outward unit normal of \f \pa\Omega_{\ell} \A at \f x $, and
     \f x \cdot n^{(\ell)} \A is the inner product of \f x \A and
     \f n^{(\ell)} \A in \f \RN $.                  
\BP

        {\bf Remark 2.2.}\ The condition (2.2) requires that the inner 
     products \f x \cdot n^{(1)} \A and \f x \cdot n^{(2)} \A do not
     change their signs almost always on \f S $. Note that the above 
     assumption is satisfied if \f \Omega_1 \A is a cylindrical domain,
     \f \mu_1 < \mu_2 $, and the origin is put in \f \Omega_1 $. Also (2.2)
     is satisfied when \f S \A is an ($ N-1$-dimensional) plane.
\BP
      
        {\bf Definition 2.3.}\ Let \f X \A be the Hilbert space defined 
     by
$$
                         X = L_2(\RN; \mu(x)dx)             \T (2.3)
$$
     with its inner product \f (\ , \ )_X \A and \f \M \ \M_X \A 
     given by
$$        
\left\{ \split
     {\dm  (f, g)_X = \int_{\RN} f(x)\overline{g(x)}\mu(x)\, dx, }\\ *[6pt]
     {\dm  \M f \M_X = [(f, f)_X]^{1/2}, }                                 
\endsplit \right.                                       \T (2.4)
$$
     Then the operator \f H \A in \f X \A is defined by (1.4), that is,
$$
\left\{ \split
            D(H) = H^2(\RN), \\  *[4pt]
            Hu = hu,
\endsplit \right.                                       \T (2.5)
$$
     where \f h \A is given by (1.1). It is easy to see that \f H \A
     is a selfadjoint operator in \f X $.

$$
\ \ \ \ \ 
$$

\cen{    {\bf \S3. The uniqueness of the solution.} }   
\BP

        In this and the following three sections we assume that \f N \ge 3$.

        In order to discuss the uniqueness of the solution of
     the inhomogeneous equation
$$
             -\, \mu(x)^{-1}\Delta u - \lambda u = f  \qq (\lambda > 0)
                                                               \T (3.1)
$$
     with radiation condition, we shall start with some notations.
\BP

        {\bf Notation 3.1.}\ \ Let \f z \in \bC $, \f x = (x_1, x_2, \cdots,
     x_N) $, \f r = |x| $, \f \tx = (\tx_1, \tx_2, \cdots, \linebreak \tx_N)$ 
     $ = x/r$, \f \pa_j = \pa/\pa x_j \A and \f \gr = (\pa/\pa x_1, 
     \pa/\pa x_2, \cdots, \pa/\pa x_N)$. Then we set
\SP

        (1) \f k = k(x) = k(x, z) = [z\mu(x)]^{1/2}$, where the branch is 
     taken so that \f {\rm Im}\,k(x, z) \ge 0 $;       
\SP

        (2) \f a = a(x) = a(x, z) = {\rm Re}\,k(x, z)$;
\SP

        (3) \f b = b(x) = b(x, z) = {\rm Im}\,k(x, z)$;                           
\SP

        (4) \f \CD_j u = \pa_j u + \{(N-1)/(2r)\}\tx_j u - ik(x)\tx_j u$,
     where \f j = 1, 2, \cdots, N$;
\SP

        (5) \f \CD u = \gr u + \{(N-1)/(2r)\}\tx u - ik(x)\tx u$;                 
\SP

        (6) \f \CD_r u = \CD u \cdot \tx = \pa u/\pa r 
                                + \{(N-1)/(2r)\} u - ik(x) u$;                 
\SP

        (7) \f \CD_n u = \CD u \cdot n = \pa u/\pa n 
                             + \{(N-1)/(2r)\} (\tx \cdot n)u 
                                 - ik(x) (\tx \cdot n)u$, 
     where \f n \A is a unit vector in \f \RN$.                                      
\BP

        Let \f u \in H^2(\RN)_{{\rm loc}} $. Then the restrictions
     \f u|_G \A and \f \pa_ju|_G $, \f j = 1, 2, \cdots, N $, of \f u \A
     and \f \pa_ju = \pa u/\pa x_j \A onto a smooth surface \f G \A 
     are defined as the traces of \f u \A and \f \pa_ju \A on \f G $, 
     respectively. Thus \f u|_G \A and \f \pa_ju|_G \A are considered to 
     belong to \f L_2(G)_{{\rm loc}} $.  

       In this section we are going to prove the following theorem:        
\BP

        {\bf Theorem 3.2.} {\it Suppose that \r{\em Assumption 2.1}$ with 
     \f N \ge 3 \A holds. Let \linebreak \f u \in H^2(\RN)_{\mb{\rm loc}} \A 
     be a solution of the homogeneous equation 
$$
           -\, \mu(x)^{-1}\Delta u - \lambda u = 0 \qq (\lambda > 0)
                         					\T (3.2)
$$
     on \f \RN \A such that                           					
$$
     \liminf_{R\to\infty} 
         \int_{S_R} \big( \bigg|\frac{\pa u}{\pa r}\bigg|^2 
                        + |u|^2 \big)\, dS = 0,                     \T (3.3)
$$                        
     where
$$     
          S_R = \{ x \in \RN \ :\ |x| = R \}.                       \T (3.4)
$$
     Then \f u \A is identically zero.}
\BP
     
        The proof will be divided into several steps. First, we are going to 
     show an identity which directly follows from the equation 
     \f -\, \mu^{-1}\Delta u - zu = f $.
\BP

        {\bf Proposition 3.3}\ ({\it cf.} \ [\IS], Lemma 2.2 {\it and}
     \ [\Sc], Lemma 2.5.) {\it Set
$$
             f = \mu(x)^{-1}(- \Delta u - k^2 u),           \T (3.5)
$$
        where \f u \in H^2(\RN)_{{\rm loc}} $. Let \f \xi \A be a 
     real-valued, continuous function on \f [0, \infty) \A \linebreak such 
     that \f \xi \A has piecewise continuous derivative. Set 
     \f \varphi(x) = \alpha(x)\xi(|x|)$, where \f \alpha \A is a simple 
     function which is constant on each \f \Omega_{\ell} $. For 
     \f 0 < r < R < \infty $, set
$$
         B_{rR} = \{ x \in \RN \ : \ r < |x| < R \ \},             \T (3.6)
$$
     Then we have     
$$
\split
    & {\dm \int_{B_{rR}} \big(b\varphi + \frac12\frac{\pa\varphi}{\pa r}\big)
                                                            |\CD u|^2\, dx 
          + \sum_{\ell=1}^{2} \int_{\pa\Omega_{\ell}\cap B_{rR}}
       \varphi {\rm Im}\big\{\overline{k}\frac{\pa u}{\pa n}
                                 \overline{u}\big\}\, dS } \\  *[6pt]  
    & \hspace{4cm}  {\dm  + \int_{B_{rR}} \big(\frac{\varphi}r 
                                - \frac{\pa\varphi}{\pa r}\big)
                               (|\CD u|^2 - |\CD_ru|^2)\, dx } \\ *[6pt]
    & \hspace{4cm}  {\dm   + c_N \int_{B_{rR}} r^{-2}\big(\frac{\varphi}r 
              - 2^{-1}\frac{\pa\varphi}{\pa r} + b\varphi\big)|u|^2\, dx  } 
                                                                   \\ *[6pt]
    & {\dm = {\rm Re} \int_{B_{rR}} \varphi\mu(x)f\overline{\CD_ru}\, dx } 
                                                                   \\ *[6pt]
    &  \ \ \ \  {\dm  + 2^{-1}\sum_{\ell=1}^2 \int_{\pa\Omega_{\ell}
                  \cap B_{rR}} \varphi \big\{ \frac{(N-1)b}r + |k|^2 \big\}
                                      (\tx \cdot n)|u|^2\, dS } \\ *[6pt]
    &  \ \ \ \  {\dm  + 2^{-1} \int_{S_R} \varphi \big(2|\CD_ru|^2 
                      - |\CD u|^2  - c_Nr^{-2}|u|^2 \big) \, dS } \\  *[6pt]
    &  \ \ \ \  {\dm - 2^{-1} \int_{S_r} \varphi(2|\CD_ru|^2 - |\CD u|^2 
                                               - c_Nr^{-2}|u|^2)\, dS, }
\endsplit                                                          \T (3.7)                                                             
$$
     where \f \Omega_1, \Omega_2 \A satisfies} \ (1.3), 
     {\it \f S \A is as in} \ (ii) {\it of} \ Assumption 2.1,
     {\it \f \pa/\pa n \A in the integrand of the surface 
     integral over \f \pa\Omega_{\ell}\cap B_{rR} \A means the directional 
     derivative in the direction of the outward normal \f n = n^{(\ell)} \A
     of \f \pa\Omega_{\ell}$, and}
$$     
              c_N = (N-1)(N-3)/4.                                \T (3.8)
$$  
\BP
      
        To prove Proposition 3.3 we first rewrite (3.5) as
$$
      - \sum_{j=1}^N \pa_j\CD_j u 
              + \big\{ \frac{N - 1}{2r} - ik \big\} \CD_r u 
                            + \frac{c_N}{r^2}u = \mu(x)f         \T (3.9)
$$
     Then (3.7) is obtained by multiplying bothe sides of (3.9) by 
     \f \varphi\overline{\CD_ru} $, taking the real part and using partial
     integration. For the details of computation see Appendix, A.1 of [\JS].
     The following lemmas are also used for the proof of Theorem 3.2.
\BP
        
        {\bf Lemma 3.4.}\ {\it Let \f \varphi(x) = \xi(|x|) \A and let 
     \f \xi \A be a continuous function on \f [0, \infty) \A such that 
     \f \xi \A has piecewise bounded continuous derivative 
     \f \xi' \A and \f \xi(0) = 0 $. Let \f S \A 
     be an $N-1$-dimensional continuous surface which consists of a finite 
     number of smooth surfaces. Let \f F(x) \A be a locally \f L_1 \A 
     function with locally \f L_1 \A derivatives in a neighborhood of \f S $.
     Then we have, \f R > 0 $,
$$
           \int_{B_R\cap S} \varphi(|x|)F(x)\, dS
                = \int_0^R \frac{\pa\varphi}{\pa r} 
                      \bigg( \int_{B_{rR}\cap S} F(x)\, dS \bigg)\, dr,
                                    			        \T (3.10)
$$
     where, for \f 0 < r < R < \infty $, \f B_{rR} \A is as in {\rm (3.6)},
     and \f B_R = \{ x \in \RN \ :\  |x| < R \} \A is an open ball with
     origin \f 0\A and radius \f R $. }       
          
\BP

        {\it Proof. \ } Since \f F(x) \A can be approximated by a 
     sequence of \f C^1 \A functions in a neighborhood of \f S $, we may 
     assume that \f F \A is a \f C^1 \A function. For \f \epsilon > 0 \A
     set 
$$
       S_{\epsilon} = \{ x \in \RN \ : \ {\rm dist}(x, S) < \epsilon \}
                                                                  \T (3.11)
$$
     where \f {\rm dist}(x, S) \A is the distance between \f x \A and
     \f S $, and let \f \chi_{S, \epsilon}(x) \A be the characteristic 
     function of the set \f S_{\epsilon} $. Then we have by an easy 
     computation
$$
\split
  & {\dm \int_0^R \frac{\pa\varphi}{\pa r} 
                   \bigg( \int_{B_{rR}\cap S_{\epsilon}} F(x)\, dx 
                                                \bigg)\, dr }  \\ *[6pt]
  & \hspace{3cm} {\dm = \int_0^R \xi'(r)                                
                 \bigg( \int_{B_{rR}} \chi_{S, \epsilon}(x)F(x)\, dS 
                                             \bigg)\, dr }  \\ *[6pt] 
  & \hspace{3cm} {\dm = \int_0^R \xi'(r) \int_r^R
                \bigg( \int_{S_t} \chi_{S, \epsilon}(x)F(x)\, dS 
                                     \bigg)\, dt\, dr } \\ *[6pt]         
  & \hspace{3cm} {\dm = \int_0^R \int_0^t \xi'(r)       
                    \bigg( \int_{S_t} \chi_{S, \epsilon}(x)F(x)\, dS 
                                           \bigg)\, drdt }  \\ *[6pt]
  & \hspace{3cm} {\dm = \int_0^R \xi(t)
                    \bigg( \int_{S_t} \chi_{S, \epsilon}(x)F(x)\, dS 
                                            \bigg)\, dt } \\ *[6pt]         
  & \hspace{3cm} {\dm = \int_{B_R} \varphi(|x|) \chi_{S, \epsilon}(x) 
                                                           F(x)\, dx. }   
\endsplit                                                          \T (3.12)
$$
     The equality (3.10) (for smooth \f F $) is obtained by dividing both 
     sides of (3.12) by \f \epsilon \A and letting \f \epsilon \downarrow 
     0 $, which completes the proof. \ \ $\Vert$
\BP

        {\bf Lemma 3.5.}\ {\it Let \f u \in H^2(\RN)_{{\rm loc}} \A 
     be a solution of the homogeneous equation} \ (3.2) {\it with 
     \f \lambda > 0 $. Let \f \Omega_1 \A and \f \Omega_2 \A satisfy} \ 
     (i) \ {\it and} \ (ii) {\it of} \ Assumption 2.1.

        (i) {\it Let \f 0 < r < R < \infty $. Then we have
$$
\split
   {\dm \sum_{\ell=1}^2 \int_{\pa\Omega_{\ell}\cap B_{rR}} 
          {\rm Im}\big\{ k\frac{\pa u}{\pa n} \overline{u} \big\}\, dS } 
                                                                 \\ *[4pt]
    \hspace{2.5cm} {\dm  = - \int_{S_R}{\rm Im}\big\{ k\frac{\pa u}{\pa r}
          \overline{u} \big\}\, dS
        + \int_{S_r} {\rm Im}\big\{ k\frac{\pa u}{\pa r}\overline{u} 
                                                          \big\}\, dS, }
\endsplit  							                                       \T (3.13) 
$$
     where \f B_{rR} \A is given by} \ (3.6).
     
        (ii) {\it Let \f \varphi(x) = \xi(|x|) \A and let \f \xi \A be a 
     real-valued, continuous function on \f [0, \infty) \A such that 
     \f \xi \A has piecewise continuous derivative \f \xi' \ge 0 $,  and 
     \f \xi(0) = 0 $. Then we have, for \f 0 < R < \infty $,
$$
     \int_{B_R} \frac12 \frac{\pa\varphi}{\pa r} |\CD u|^2\, dx  
                + \sum_{\ell=1}^2 \int_{\pa\Omega_{\ell}\cap B_R}
                    \varphi  {\rm Im}\big\{k\frac{\pa u}{\pa n}
                                 \overline{u}\big\}\, dS   \hspace{3cm}
$$
$$ 
         \ge \int_{B_R} \frac12 \frac{\pa\varphi}{\pa r} k^2 |u|^2\, dx 
              + \int_{B_R} \frac12 \frac{\pa\varphi}{\pa r} ( |\gr u|^2 
                                - \bigg|\frac{\pa u}{\pa r}\bigg|^2)\, dx
$$
$$
              - \xi(R) \int_{S_R} {\rm Im}\big\{ k\frac{\pa u}{\pa r}
                                 \overline{u} \big\}\, dS,        \T (3.14)
$$
     where \f k = k(x) = \sqrt{\lambda\mu(x)} $.}  
\BP

        {\it Proof. \ } (I) Multiply both side of 
$$
          -\, \Delta u - \mu(x)\lambda u = 0                     \T (3.15)
$$            
     by \f k\overline{u} \A and integrate over \f B_{rR} $. Then, (3.13) is 
     obtained by taking the imaginary part and using partial integration.
\SP
               					
        (II) It follows from (3.13) and Lemma 3.4 that
$$
  \sum_{\ell=1}^2 \int_{\pa\Omega_{\ell}\cap B_R} \varphi
         {\rm Im}\big\{ k\frac{\pa u}{\pa n} \overline{u} \big\}\, dS 
                                                          \hspace{5.5cm}
$$
$$
         = \sum_{\ell=1}^2 \int_0^R \frac{\pa\varphi}{\pa r}  
            	\bigg( \int_{\pa\Omega_{\ell}\cap B_{rR}}									
	{\rm Im}\big\{ k\frac{\pa u}{\pa n} \overline{u} \big\}\, dS \bigg)
	                                              	\, dr  \hspace{2cm}					
$$
$$
       = - \int_0^R \frac{\pa\varphi}{\pa r} \bigg(		                                          					
	   \int_{S_R} {\rm Im}\big\{ k\frac{\pa u}{\pa r}\overline{u} 
	     \big\}\, dS \bigg)\, dr                            \hspace{3cm}	
$$
$$
         	   + \int_0^R \frac{\pa\varphi}{\pa r} \bigg(		                                          					
	   \int_{S_r} {\rm Im}\big\{ k\frac{\pa u}{\pa r}\overline{u} 
	     \big\}\, dS \bigg)\, dr                                  
$$
$$
      = -\,\xi(R) \int_{S_R} {\rm Im}\big\{ k\frac{\pa u}{\pa r}
               \overline{u} \big\}\, dS	     
  	+ \int_{B_R} \frac{\pa\varphi}{\pa r}
  	       {\rm Im}\big\{ k\frac{\pa u}{\pa r}\overline{u} \big\}\, dx. 
  	                                                       \T (3.16)					
$$

        (III) Now we are going to evaluate the term \f |\CD u|^2 $. By 
     definition and the Schwarz inequality it follows that
$$        										
\split
       {\dm |\CD u|^2 = |\gr u|^2 - 2{\rm Im}\big( k\frac{\pa u}{\pa r}
                                  \overline{u} \big) + k^2|u|^2 } \\ *[5pt]    
       \hspace{3.5cm} {\dm + \frac{N-1}r{\rm Re}\big( \frac{\pa u}{\pa r}
               \overline{u} \big) + \frac{(N-1)^2}{4r^2}|u|^2. }   \\ *[5pt] 
       \hspace{1cm}  {\dm \ge - 2{\rm Im}\big( k\frac{\pa u}{\pa r} 
                                   \overline{u} \big) 
                     + ( |\gr u|^2 - |\frac{\pa u}{\pa r}|^2) + k^2|u|^2. }  
\endsplit                                                         \T (3.17)
$$
     Multiply both side of (3.17) by \f \pa\varphi/\pa r \A and integrate 
     over \f B_R $. Then, using (3.16), too, we have (3.14).  \ \ $\Vert$                          
\BP
        
        Now we are in a position to prove Theorem 3.2. Here and in the sequel
     we agree that \f C = C(A, B, \cdots) \A in an inequality means a
     positive constant depending on \f A, B, \cdots $. But very often 
     symbols indicating obvious dependence such as the operator \f H \A will
     be left out. 
\BP

        {\it Proof of} \ Theorem 3.2. \ Let 
     \f u \in H^2(\RN)_{{\rm loc}} \A be a solution of the homogeneous
     equation (3.2). Let \f R_0 \ge 1 \A and define \f \varphi \A by
$$
       \varphi(x) =
\left\{ \split
           |x|  \q (0 \le |x| \le R_0),  \\ *[4pt]
           R_0  \q (|x| > R_0).
\endsplit \right.                                           \T (3.18)
$$                 
     The function \f \varphi \A satisfies the conditions given in
     Proposition 3.3 and Lemma 3.5. Then it follows from (3.7) in 
     Proposition 3.3 with \f f = 0 $, \f b = 0 $, \f k = 
     \sqrt{\lambda\mu(x)} \A and Lemma 3.5, (ii) that, for any 
     \f R > R_0 > r > 0 $, 
$$
\split
      {\dm \int_{B_{rR_0}} \frac12 k^2|u|^2\, dx } \\  *[6pt] 
      \hspace{1cm} {\dm \le R_0 C\int_{S_R} \big( \big|\frac{\pa u}{\pa r}
                          \big|^2 + |u|^2 \big)\, dS } \\  *[6pt]            
      \hspace{1.5cm} {\dm + 2^{-1} \int_{S_r} r( \big|\CD u\big|^2 
                                + c_Nr^{-2}|u|^2)\, dS  } \\  *[6pt]  
      \hspace{1.5cm} {\dm + \int_{B_r} \frac12 r|\CD u|^2\, dx
                    + \sum_{\ell=1}^{2} \int_{\pa\Omega_{\ell}\cap B_r}
                         r{\rm Im}\big\{ k\frac{\pa u}{\pa n}
                                 \overline{u}\big\}\, dS }      
\endsplit                                                         \T (3.19)
$$
      with a positive constant \f C = C(\lambda) $, where we have used the 
      facts that
$$
\left\{ \split
     {\dm \frac{\varphi}r - 2^{-1}\frac12\frac{\pa\varphi}{\pa r} \ge 0, } \\  *[5pt]
     {\dm \big|\gr u\big|^2 - \bigg|\frac{\pa u}{\pa r}\bigg|^2 \ge 0, } 
                                                                 \\ *[5pt]
     {\dm 2|\CD_ru|^2 - |\CD u|^2 - c_Nr^{-2}|u|^2 } \\ *[5pt]
     {\dm \ \ \ \ \ \  = |\CD_ru|^2 - (|\CD u|^2 - |\CD_ru|^2) 
                          - c_Nr^{-2}|u|^2 \le |\CD_ru|^2,  } \\ *[5pt]
     {\dm - 2|\CD_ru|^2 + |\CD u|^2 + c_Nr^{-2}|u|^2 
                               \le |\CD u|^2 + c_Nr^{-2}|u|^2, }  \\ *[5pt]
     {\dm \sum_{\ell=1}^2 \int_{\pa\Omega_{\ell}\cap B_R}
                     \varphi |k|^2(\tx \cdot n)|u|^2\, dS  }  \\  *[5pt] 
     {\dm \ \ \ \ \ \                 
         = \lambda \int_{S\cap B_R} \varphi(\mu_1 - \mu_2)
                           (\tx\cdot n^{(1)})|u|^2 \, dS \le 0  }    
\endsplit \right.                                                 \T (3.20)
$$
     Here the last inequality follows from (2.2) in Assumption 2.1. 
     Since \f u \in H^2(\RN)_{{\rm loc}} \A with \f N \ge 3 $, it follows 
     from the Hardy inequality that \f u/r \A is locally \f L_2(\RN) $, 
     and hence the second term of the right-hand side of (3.19) tends to 
     zero as \f r \to 0 \A along a suitable sequence \f \{ r_m \} $, i.e., 
$$
     \int_{S_{r_m}} r( \big|\CD u\big|^2 + c_Nr^{-2}|u|^2)\, dS
                                                           \to 0.   \T (3.21)                        
$$
     as \f r_m \to 0 $.  The last two terms of the right-hand side of (3.19)
     tend to \f 0 \A as \f r \downarrow 0$, since their integrands are 
     integrable. Thus we have
$$
    \int_{B_{R_0}} \frac12 k^2|u|^2\, dx 
         \le R_0 C\int_{S_R} \big( \big|\frac{\pa u}{\pa r}\big|^2 
                                + |u|^2 \big)\, dS   \qq (R > R_0). 
                                                              \T (3.22)
$$
     Therefore, by letting \f R \to \infty \A along an appropriate sequence 
     \f \{ R_m \} $, the right-hand side of (3.22) becomes \f 0 $, i.e., 
     we have,  
$$
            \int_{B_{R_0}} \frac12 k^2|u|^2\, dx = 0            \T (3.23)
$$
     for any \f R_0 \ge 1 $, which implies that \f u \A is identically zero.
     \ \ $\Vert$            
\BP
 
        Using Theorem 3.2, we can easily show the nonexistence of the 
     eigenvalues of the operator \f H $.
\BP

        {\bf Corollary 3.6.} {\it Suppose} \ Assumption 2.1 {\it with 
     \f N \ge 3 \A holds. Then the operator \f H \A has no eigenvalues.}
\BP

        {\it Proof. \ } Since \f H \A is nonnegative, we have only
     to show that \f H \A has no nonnegative eigenvalues. Suppose that
     \f u \in D(H) = H^2(\RN) \A be an eigenfunction associated with 
     a positive eigenvalue \f \lambda \A of \f H $. Then, since \f u \A
     is a solution of (3.2) and satisfies the condition (3.3), \f u \A is
     identically zero, which is a contradiction. Suppose that \f 0 \A is an
     eigenvalue of \f H \A and \f u \A is the corresponding eigenfunction.
     Then \f u \A becomes also an eigenfunction of the operator 
     \f - \, \Delta \A and \f \lambda = 0 \A becomes an eigenvalue of
     \f - \, \Delta $, which is again a contradiction since 
     \f - \, \Delta \A does not have eigenvalue \f \lambda = 0 $. \ \ $\Vert$
\BP

        Finally we shall show that the radiation condition
$$
      \liminf_{R\to\infty} \int_{S_R} |\CD_r^{(\pm)} u|^2\, dS = 0    
                                                                   \T (3.24)
$$
     or
$$
     \liminf_{R\to\infty} \int_{S_R} |\frac{\pa u}{\pa r} 
                                        \mp iku|^2\, dS = 0        \T (3.25)
$$
     implies that (3.3) for a solution \f u \in H^2(\RN)_{{\rm loc}} \A 
     of the equation (3.2). Here \f k = \sqrt{\lambda\mu(x)} \A and
     \f \CD_r^{(\pm)} u \A is given by 
$$
    \CD_r^{(\pm)} u = \pa u/\pa r + \{(N-1)/(2r)\} u \mp ik(x) u.
                                                                  \T (3.26)
$$          
\BP
        {\bf Theorem 3.7.}\ {\it Suppose} \ Assumption 2.1 {\it with 
     \f N \ge 3 \A holds. Let \f u \in H^2(\RN)_{{\rm loc}} \A 
     be a solution of the homogeneous equation} \ (3.2) {\it with 
    \f \lambda > 0 $. Suppose that} \ (3.24) {\it or} \ (3.25) {\it holds. 
     Then \f u \A is identically zero.}
\BP

        {\it Proof. \ } We have only to show that the condition (3.24) or 
     (3.25) implies (3.3). Here we shall consider the condition
$$
      \liminf_{R\to\infty} \int_{S_R} |\CD_r^{(+)} u|^2\, dS = 0.    
                                                                  \T (3.27)
$$
     All the other conditions can be treated similarly. Multiply both 
     sides of (3.2) by \f \overline{u} $, integrate over 
     \f B_R \A with \f R > 0 \A and take the imaginary part. Then we obtain
$$
         {\rm Im} \int_{S_R} \frac{\pa u}{\pa r}\overline{u} \, dS = 0.
                                                            \T (3.28)
$$                                                            
     Since we have from (3.28)
$$
         {\rm Im} \int_{S_R} (\CD_r^{(+)} u)\overline{u} \, dS
                           = - \int_{S_R} k|u|^2 \, dS      \T (3.29)
$$
     with \f k = \sqrt{\lambda\mu} $, it follows that
$$
   k_0\int_{S_R} |u|^2 \, dS \le \int_{S_R} |\CD_r^{(+)} u||u| \, dS 
                                                       \hspace{3cm}
$$
$$
    \hspace{5cm}  \le \frac{k_0}2 \int_{S_R} |u|^2 \, dS
                         + \frac1{2k_0}\int_{S_R} |\CD_r^{(+)} u|^2 \, dS,
                                                              \T (3.30)
$$
     or
$$
     \int_{S_R} |u|^2 \, dS
           \le \frac1{k_0^2}\int_{S_R} |\CD_r^{(+)} u|^2 \, dS,
                                                         \T (3.31)
$$
     where \f k_0 = \sqrt{\lambda\min(\mu_1, \mu_2)} $. On the other hand
     we have
$$
     \int_{S_R} \bigg|\frac{\pa u}{\pa r}\bigg|^2 \, dS
         \le 2 \int_{S_R} |\CD_r^{(+)} u|^2 \, dS 
          + 2 \bigg( \big(\frac{N-1}2\big)^2 + k_1^2\bigg) \int_{S_R} 
                                                            |u|^2 \, dS
                                                          \T (3.32)
$$
     for \f R \ge 1 \A with \f k_1 = \sqrt{\lambda\max(\mu_1, \mu_2)} $. 
     Thus it follows from (3.31) and (3.32) that there exists a positive
     constants \f C = C(\lambda) \A such that                                                          
$$
   \int_{S_R} \big( \bigg|\frac{\pa u}{\pa r}\bigg|^2 + |u|^2 \big)\, dS       
                   \le C \int_{S_R} |\CD_r^{(+)} u|^2 \, dS,    \T (3.33)
$$
     which completes the proof. \ \ $\Vert$
\BP

       Later we shall need the following corollary which guarantees the 
     uniqueness of the inhomogenous equation
$$
         - \, \mu(x)^{-1}\Delta u - \lambda u = f                  \T (3.34)
$$
     with one of the conditions
$$
\left\{ \split                                                       
     {\dm \int_{E_R} \frac1r |\CD_r^{(\pm)} u|^2 \, dx < \infty, }  
                                                                \\  *[6pt]
     {\dm \int_{E_R} \frac1r |\frac{\pa u}{\pa r} \mp iku|^2, dx < \infty, }
\endsplit \right.                                                  \T (3.35)       
$$       
     where \f E_R = \{ x \in \RN \ : \ |x| > R \}, \ R > 0 $.
\BP

        {\bf Corollary 3.8.} {\it Let \f \lambda > 0 \A and let 
     \f f \in L_2(\RN)_{{\rm  loc}} $. Then the solution 
     \f u \in H^2(\RN)_{{\rm  loc}} \A of the equation} \ (3.34) {\it with 
     one of the radiation conditions in} \ (3.35) {\it is unique.}
\BP

        {\it Proof. \ } Let \f u_1 \A and \f u_2 \A be the solutions of
     the equation (3.34) satisfying, say,
$$        
        \int_{E_R} \frac1r |\CD_r^{(+)} u_j|^2 \, dx < \infty  
                          \qq (j = 1, 2)                         \T (3.36)
$$
     with \f R > 0 $. Set \f u = u_1 - u_2. $. Then \f u \A is a solution of 
     the homogeneous equation (3.2) and satisfies (3.36) with \f u \A 
     replaced by \f u_1 - u_2 $, which implies that                              
$$
      \liminf_{R\to\infty} \int_{S_R} |\CD_r^{(+)} u|^2\, dS = 0.   
                                                                   \T (3.37)
$$
     Therefore we have \f u = 0 \A from Theorem 3.7. \ \ $\Vert$

$$
\ \ \ \ \
$$

\cen{      {\bf \S4. An estimate for the radiation condition.} }   
\SP

        Let \f L_{2, t}(\RN) \A be the weighted Hilbert space defined by 
     (1.11). Let the resolvent \f (H - z)^{-1} \A of the operator \f H \A 
     will be denoted by \f R(z) $. Now consider \f u \in X \A defined by    
$$
\left\{ \split
          u = R(z)f,  \\ *[4pt]
          z = \lambda + i\eta  \qq (\lambda \ge 0, \eta \ne 0),  \\ *[4pt]
          f \in L_{2, \delta}(\RN).
\endsplit \right.                                            \T (4.1)
$$
     In this section we are going to prove the following
\BP

        {\bf Theorem 4.1.} {\it Suppose that \r{\em Assumption 2.1}$ 
     holds with \f N \ge 3 $. Let \f 1/2 < \delta \le 1 $. Let \f u \A be 
     given by \r{\em (4.1)}$. Then there exists a positive constant 
     \f C = C(\delta) \A such that
$$
        \M \CD u \M_{\delta-1} \le C \M f \M_{\delta},
                                                                \T (4.2)   
$$
     where \f \CD u \A is as in \r{\em Notation 3.1}$, \f \M \ \M_t \A is 
     the norm of \f L_{2, t}(\RN) $, and the constant \f C(\delta) \A is 
     independent of \f f \A and \f z \A satisfying \r{\em (4.1)}$.}
\BP

        In order to show the theorem we need a lemma.
\BP

        {\bf Lemma 4.2.} {\it Let \f z = \lambda + i\eta \in \bC \bss \Ro $. 
     Let \f a $, \f b $, \f k \A be as in \r{\em Notation 3.1.\ }$ Then 
     we have}
$$
\left\{ \split
    {\dm a = a(x, z) = \sqrt{\mu(x)}c_a(z)  \q
      \bigg(c_a(z) = \frac{\eta}{|\eta|}\sqrt{\frac{|z| + \lambda}2}\bigg),} 
                                                                    \\ *[4pt]
    {\dm b = b(x, z) = \sqrt{\mu(x)}c_b(z) \q
           \bigg(c_b(z) = \frac{|\eta|}{\sqrt{2(|z| + \lambda)}}\bigg), } 
                                                                   \\ *[4pt]
        {\dm |k|^2 = |k(x, z)|^2 = \mu(x)|z|. }   
\endsplit \right.                                                               \T (4.3)
$$
\BP

        Since the lemma is shown by an easy computation, we shall omit
     the proof.
\BP

        {\it Proof of \r{\em Theorem 4.1.\ }$} (I) Let \f \xi(r) \A
     be defined by
$$
     \xi(r) =
\left\{ \split
         r           \qq \q \hspace{1.5cm} (0 \le r \le 1), \\ *[4pt]
         2^{-(2\delta-1)}(1 + r)^{2\delta-1} \q (r \ge 1)
\endsplit \right.                                             \T (4.4)                  
$$
     Set \f \varphi(x) = \xi(r)/\sqrt{\mu(x)} $, where \f r = |x| $, and 
     \f u = R(z)f \A in (3.7) of Proposition 3.3. We are 
     going to evaluate each term of the left-hand side and the right-hand 
     side of (3.7). Here we set \f 0 < r < 1 < R $.
\SP

        (II) Let
$$
      I_{L1}
        = \int_{B_{rR}} \big(b\varphi + \frac12\frac{\pa\varphi}{\pa r}\big)
                                           |\CD u|^2\, dx.       \T (4.5)
$$
     Then we have
$$
\split
     {\dm I_{L1} \ge \frac12 \int_{B_{r1}} \frac1{\sqrt{\mu(x)}}
                   |\CD u|^2\, dx } \\ *[6pt]   
     \hspace{2.5cm} {\dm + \frac12 \int_{B_{1R}} 
                       \frac{2\delta-1}{2^{2\delta-1}\sqrt{\mu(x)}}
                        (1 + r)^{2\delta-2}|\CD u|^2\, dx  } \\ *[6pt]                                   
     \hspace{.6cm} {\dm \ge \frac{c_{\delta}}{2\sqrt{M_0}} 
                     \int_{B_{rR}} (1 + r)^{2\delta-2} |\CD u|^2\, dx, }
\endsplit                                                     \T (4.6)
$$
     where \f M_0 = \max(\mu_1, \mu_2) $, and
$$
        c_{\delta} = \frac{2\delta-1}{2^{2\delta-1}}.           \T (4.7)
$$                                                                          
\SP
                                                
        (III) Let the second term of the left-hand side of (3.7) be denoted
     by \f I_{L2} $. Note that we have   
$$
\split
     {\dm \int_{\pa\Omega_{\ell}\cap B_{rR}} 
       \varphi(x) {\rm Im}\big\{\overline{k}\frac{\pa u}{\pa n}
                                 \overline{u}\big\}\, dS  } \\ *[6pt]        
     \hspace{2cm}{\dm  = \int_{S\cap B_{rR}} \frac{\xi(r)}{\sqrt{\mu_{\ell}}} 
                            \sqrt{\mu_{\ell}}  
                  {\rm Im}\big\{c_-(z)\frac{\pa u}{\pa n^{(\ell)}}
                                 \overline{u}\big\}\, dS   } \\ *[6pt]               
     \hspace{2cm}{\dm   = \int_{S\cap B_{rR}} \xi(r) {\rm Im}\big\{c_-(z)
                 \frac{\pa u}{\pa n^{(\ell)}}\overline{u}\big\}\, dS 
                                                    \q (\ell = 1, 2),}                         
\endsplit                                                         \T (4.8)
$$
     where \f c_-(z) = c_a(z) - ic_b(z) $, and \f c_a(z) \A and 
     \f c_b(z) \A are as in Lemma 4.2. Noting that 
     \f n^{(1)} +  n^{(2)} = 0 $, we have \f I_{L2} = 0 $.
\SP

        (IV) Set
$$
 I_{L3} = \int_{B_{rR}} \big(\frac{\varphi}r - \frac{\pa\varphi}{\pa r}\big)
               (|\CD u|^2 - |\CD_ru|^2)\, dx.                    \T (4.9)         
$$
     Since it is easy to see that the integrand is nonnegative, we have     
     \f I_{L3} \ge 0 $. Similarly we have
$$
    I_{L4} = c_N \int_{B_{rR}} r^{-2}\big(\frac{\varphi}r  
               - 2^{-1}\frac{\pa\varphi}{\pa r} + b\varphi\big)|u|^2\, dx
               \ge 0.                                             \T (4.10)
$$
\SP

        (V) Using the Schwarz inequality, we have
$$
\split
     {\dm I_{R1} = {\rm Re} \int_{B_{rR}} 
                   \varphi\mu(x)f\overline{\CD_ru}\, dx } \\ *[6pt] 
     \hspace{1cm} {\dm \le \sqrt{M_0} \int_{B_{rR}} 
                     (1 + r)^{2\delta-1}|f|\CD_ru|\, dx } \\ *[6pt] 
     \hspace{1cm} {\dm  \le \frac{M_0}{4\epsilon} 
            \int_{B_{rR}} (1 + r)^{2\delta}|f|^2 \, dx        
              + \epsilon \int_{B_{rR}} 
                   (1 + r)^{2\delta-2}|\CD u|^2 \, dx, }
\endsplit                                                       \T (4.11)
$$
     where \f \epsilon \A is an arbitrary positive number.                    
\SP

        (VI) Let \f I_{R2} \A be the second term of the right-hand side of
     (3.7). Then,
$$
\split
     {\dm I_{R2} = 2^{-1} \bigg( \sum_{\ell=1}^2 \int_{\pa\Omega_{\ell}
                                                              \cap B_{rR}}
                  \varphi \frac{(N-1)b}r (\tx \cdot n)|u|^2\, dS } \\ *[6pt]
     \hspace{3cm} {\dm + \sum_{\ell=1}^2 \int_{\pa\Omega_{\ell}\cap B_{rR}}
                   \varphi|k|^2 (\tx \cdot n)|u|^2\, dS \bigg)  } \\ *[6pt]
     \hspace{.7cm} {\dm \equiv 2^{-1} \big[ I_{R2}^{(1)} + I_{R2}^{(2)} 
                                                                   \big]. }
\endsplit                                                         \T (4.12)
$$
     Here, as in (II), we see from Lemma 4.2 that \f I_{R2}^{(1)} = 0 \A and

$$
      I_{R2}^{(2)}      
           = \int_{S\cap B_{rR}} \frac{\xi(r)}{\sqrt{\mu_1}} 
                  \mu_1|z| (\tx \cdot n^{(1)}) |u|^2\, dS   \hspace{3.2cm}
$$
$$
            \ \ \ \ \ \ \ \ \ \ \ \ \ \ \ \
                     + \int_S \frac{\xi(r)}{\sqrt{\mu_2}}
                     \mu_2|z| (\tx \cdot n^{(2)}) |u|^2\, dS  \hspace{2.1cm}
$$
$$
            = \int_{S\cap B_{rR}} \xi(r) \big( \sqrt{\mu_1} - 
                   \sqrt{\mu_2} \big) |z| (\tx \cdot n^{(1)}) |u|^2\, dS 
                                                                \hspace{.7cm}
$$
$$
            = \int_{S\cap B_{rR}} \frac{\xi(r)}{\sqrt{\mu_1}+\sqrt{\mu_2}}
                     (\mu_1 - \mu_2) |z| (\tx \cdot n^{(1)}) |u|^2\, dS 
$$
$$
            \le 0,    \hspace{7.7cm}
                                                        \T (4.13)
$$
     where we have used (2.2) in Assumption 2.1 again.  Thus we have
     \f I_{R2} \le 0 $.
\SP

        (VII) It follows from (II) $\sim$ (VI) that   
        
\newpage
                                  
$$
    \frac{c_{\delta}}{4\sqrt{M_0}} 
               \int_{B_{rR}} (1 + r)^{2\delta-2} |\CD u|^2\, dx \hspace{5cm}
$$
$$
        \le \frac{M_0^{3/2}}{c_{\delta}} 
          \int_{B_{rR}} (1 + r)^{2\delta}|f|^2 \, dx  \hspace{2.5cm}
$$
$$
         \ \ \ \ \ \ \ + 2^{-1} \int_{S_R} \varphi 
                             ( |\CD_ru|^2 \, dS   \hspace{2.8cm}  
$$
$$
         \ \ \ \ \ \ \ \ \ \ \ \
                          + 2^{-1} \int_{S_r} \varphi(|\CD u|^2
                                   + c_Nr^{-2}|u|^2) \, dS,           
                                                         \T (4.14)
$$                  
     where we set 
$$
           \epsilon = \frac{c_{\delta}}{4\sqrt{M_0}}            \T (4.15)
$$
     in (4.11), and we have used the third and fourth inequalities in (3.20) 
     to eveluate the third and fourth term of the right-hand side of (3.7), 
     Proceeding as in the proof of Theorem 3.2, we can let 
     \f R \uparrow \infty \A and \f r \downarrow 0 \A along suitable 
     sequences \f \{ R_n \} \A and \f \{ r_n \} $ to obtain
$$
     \frac{c_{\delta}}{2\sqrt{M_0}} 
              \int_{\RN} (1 + r)^{2\delta-2} |\CD u|^2\, dx  \hspace{3.5cm}
$$
$$
        \le \frac{M_0^{3/2}}{c_{\delta}} 
            \int_{\RN} (1 + r)^{2\delta}|f|^2 \, dx,
                                                                  \T (4.16)
$$                  
     which completes the proof. \ \ $\Vert$
\BP

        Theorem 4.1 is combined with the inequality
$$
        |\gr u - ik\tx u|^2 \le 2|\CD u|^2 + \frac{(N-1)^2}{2r^2}|u|^2
                                                                  \T (4.17)
$$
     to obtain the following corollary:
\BP

        {\bf Corollary 4.3.} {\it Suppose that \r{\em Assumption 2.1}$ with 
     \f N \ge 3 \A holds. Let \f 1/2 < \delta \le 1 $. Let \f u \A be given 
     by \r{\em (4.1)}$. Then there exists a positive constant 
     \f C = C(\delta) \A such that
$$
        \int_{E_1} (1 + r)^{2\delta-2} |\gr u - ik\tx u|^2 \, dx
                \le C \big( \M f \M_{\delta}^2 + \M u \M_{-\delta}^2 \big),
                                                                  \T (4.18)
$$
     where} 
$$
            E_1 = \{ x \in \RN \ : \ |x| > 1 \ \}.                \T (4.19)
$$ 

$$
\ \ \ \
$$

\cen{     {\bf \S5. Boundedness of \f R(z) $.}}        
\BP

        Using the estimates for the radiation condition term \f \CD u \A
     (\f u = R(z)f $), which were given in the preceding section, we are 
     going to prove several uniform boundedness estimates for \f R(z) $.
     At the same time the first theorems (Theorems 5.1 and 5.2) will  
     prepare the arguments given in \S6, where we shall discuss the
     limiting absorption principle for the operator \f H $.  
\BP

        {\bf Theorem 5.1.} {\it Suppose that \r{\em Assumption 2.1}$ holds
     with \f N \ge 3 $. Let \f 1/2 < \delta \le 1 $. Let \f u = R(z)f $, 
     where \f f \in L_{2,\delta}(\RN) \A and \f z = \lambda + i\eta \A with
     \f \lambda \ge 0 $, \f \eta \ne 0 $. Then there exists a positive 
     constant \f C = C(\delta) \A such that 
$$
    \int_{E_s} (1 + r)^{-2\delta} 
                    \big( |\gr u|^2 + |k|^2|u|^2 \big) \, dx  \hspace{6cm}
$$
$$
          \le C(1 + \sqrt{|z|})(1 + s)^{-(2\delta-1)} 
                    \big( \M f \M_{\delta}^2 + \M u \M_{-\delta}^2 \big)
                                 \q (s \ge 1),  
                                                      \T (5.1)
$$
     where \f E_s = \{ x \in \RN \ : \ |x| > s \} $, and \f \M \ \M_{t} \A 
     is the norm of \f L_{2,t}(\RN) $. The constant \f C = C(\delta) \A is 
     independent of \f f \A and \f z \A satisfying the above conditions 
     and \f s \ge 1 $.}                                                 
\BP

        {\it Proof. \ } Let \f \alpha(x) = 1/\sqrt{\mu(x)} $. Then. as we 
     have seen in  (III) of the proof of Theorem 4.1, we have
$$
        \sum_{\ell=1}^{2} \int_{\pa\Omega_{\ell}\cap B_{r}}
          \alpha {\rm Im}\big\{\overline{k}\frac{\pa u}{\pa n}
                                 \overline{u}\big\}\, dS = 0.     \T (5.2)
$$
     Multiply both sides of  the equation \f \mu f = - \Delta u - k^2u \A 
     by \f \alpha\overline{k}\overline{u} $, integrate over \f B_r $ and
     take the imaginary part to obtain
$$
\split
      {\dm \int_{B_r} \alpha \mu {\rm Im}
                 \big( \overline{k}f\overline{u} \big) \, dx 
                   = - \int_{B_r} b\alpha \big( |\gr u|^2 
                                  + |k|^2|u|^2 \big) \, dx } \\ *[6pt]
      \hspace{4.5cm} {\dm - \int_{S_r} \alpha{\rm Im}\big( \overline{k}
                             \frac{\pa u}{\pa r}\overline{u} \big) \, dS, }
\endsplit                                                         \T (5.3)
$$
     where we have used (5.2). Combining
$$
  |\gr u - ik\tx u|^2 
        =  |\gr u|^2 + |k|^2|u|^2 
               - 2{\rm Im}\big( \overline{k}
                          \frac{\pa u}{\pa r}\overline{u} \big)   \T (5.4)
$$
     with (5.3), we obtain
$$
\split
    {\dm \int_{S_r} } \alpha  \big( |\gr u|^2 + |k|^2|u|^2 \big) \, dS 
                                                               \\ *[6pt]
     \hspace{1.5cm} {\dm = \int_{S_r} \alpha |\gr u - ik\tx u|^2 \, dS
                   - 2 \int_{B_r} \alpha\mu {\rm Im}
                        \big( \overline{k}f\overline{u} \big) \, dx } 
                                                                \\ *[6pt]
      \hspace{2.5cm} {\dm - \int_{B_r} b\alpha \big( |\gr u|^2 
                                 + |k|^2|u|^2 \big) \, dx }  \\ *[6pt]
      \hspace{1.5cm} {\dm \le \int_{S_r} \alpha |\gr u - ik\tx u|^2 \, dS
               + 2\sqrt{M_0} \M f \M_{\delta}\M ku \M_{-\delta}, }             
\endsplit                                                        \T (5.5)
$$
     Multiply both sides of (5.5) by \f (1 + r)^{-2\delta} \A and integrate
     from \f s \A to \f \infty $. Then, setting  
     \f \mu_0 = \min(\mu_1, \mu_2) \A and \f M_0 = \max(\mu_1, \mu_2) $, 
     we have
$$
\split
      {\dm \frac1{\sqrt{M_0}} } & {\dm \int_{E_s} (1 + r)^{-2\delta} 
                    \big( |\gr u|^2 + |k|^2|u|^2 \big) \, dx }  \\
         & {\dm \le \frac1{\sqrt{\mu_0}}\int_{E_s} (1 + r)^{-2\delta} 
                                      |\gr u - ik\tx u|^2 \, dx } \\
         & \ \ \ \ \ \ \ \ \ \ \ \ \ \ \ \ \ \ \ \ \ 
              {\dm + \frac{2\sqrt{M_0}}{(2\delta-1)}(1 + s)^{-(2\delta-1)}
                                 \M f \M_{\delta}\M ku \M_{-\delta} } \\ 
         & {\dm \le \frac{(1 + s)^{-(4\delta-2)}}{\sqrt{\mu_0}}
                \int_{E_s} (1 + r)^{2\delta-2} |\gr u - ik\tx u|^2 \, dx } \\  
         & \ \ \ \ \ \ \ \ \ \ \ \ \ \ \ \ \ \ \ \ \ 
               {\dm + \frac{2\sqrt{M_0}}{(2\delta-1)}(1 + s)^{-(2\delta-1)}
                          \M f \M_{\delta}\M ku \M_{-\delta}, } 
\endsplit                                                          \T (5.6)
$$
     which, together with Corollary 4.3, gives (5.1). \ \ $\Vert$ 
\BP

        In the next theorem an improved estimate for \f u = R(z)f \A
     will be given.
\BP

        {\bf Theorem 5.2.} {\it Suppose that \r{\em Assumption 2.1}$ holds
     with \f N \ge 3$. Let \f 1/2 < \delta \le 1 \A. Let \f u = R(z)f $, 
     where \f f \in L_{2,\delta}(\RN) \A and \f z = \lambda + i\eta \A with
     \f \lambda \ge 0 $, \f \eta \ne 0 $. Then there exists a positive 
     constant \f C = C(\delta) \A such that        
$$
    \int_{E_s} (1 + r)^{-2\delta} |u|^2 \, dx 
           \le C(1 + s)^{-2(2\delta-1)}\frac1{|z|}\M f \M_{\delta}^2  
                                           \q (s \ge 0).         \T (5.7)
$$
     The constant \f C = C(\delta) \A is independent of \f f \A and 
     \f z \A satisfying the above conditions and \f s \ge 0 $.}        
\BP

        {\it Proof. \ } Multiply both sides of  the equation 
     \f \mu f = - \Delta u - k^2u \A  by \f a\alpha\overline{u} $, 
     integrate over \f B_r $ and take the imaginary part, where
     \f \alpha(x) = 1/\sqrt{\mu(x)} \A again. Then we have
$$
\split
      {\dm \int_{B_r} a \alpha \mu {\rm Im}
                              \big(f\overline{u} \big) \, dx 
            = - \int_{S_r} a\alpha{\rm Im}\big(
                     \frac{\pa u}{\pa r}\overline{u} \big) \, dS, } \\ *[6pt]
      \hspace{5cm} {\dm - \int_{B_r} 2a^2b\alpha |u|^2 \, dx } \\ *[6pt]
      \hspace{3.2cm} {\dm \le -\int_{S_r} a\alpha {\rm Im} 
                                \big(\frac{\pa u}{\pa r} 
                                       \overline{u}\big) \, dS, }
\endsplit                                                         \T (5.8)  
$$
     where we have noted that (5.2) holds with \f \overline{k} \A replaced by \f a $.
     Since 
$$
\split
      |\CD_ru |^2 &
           = {\dm \big| \frac{\pa u}{\pa r} + \frac{N-1}{2r} + bu \big|^2
              + a^2|u|^2 
              - 2a {\rm Im}\big(\frac{\pa u}{\pa r}\overline{u}\big)  } \\         
         & {\dm \ge a^2|u|^2 
              - 2a {\rm Im}\big(\frac{\pa u}{\pa r}\overline{u}\big), }
\endsplit                                                         \T (5.9)
$$
     we have from (5.8)     
$$
    \int_{S_r} a^2\alpha|u|^2 \, dS   
             \le \int_{S_r} \alpha |\CD_ru |^2 \, dS
                  - 2 \, \int_{B_r} a\alpha \mu{\rm Im} 
                          \big(f \overline{u}\big) \, dx,          \T (5.10)
$$
     and hence, by the use of the first relation of (4.3), it follows that
$$
  \int_{S_r} |u|^2 \, dS
        \le \frac4{e(z)^2\mu_0}\int_{S_r} |\CD_ru |^2 \, dS
          + \frac{4{M_0}}{e(z)\sqrt{\mu_0}}\M f \M_{\delta}\M u \M_{-\delta}
                                                                    \T (5.11)
$$
     with \f e(z) = \sqrt{2(|z|+\lambda)} $. Multiply both sides of (5.11) 
     by \f (1 + r)^{-2\delta} \A and integrate on 
     \f (s, \infty) $. Then we see that, for \f s \ge 0 $,
$$     
\split
    {\dm \int_{E_s} } & {\dm (1 + r)^{-2\delta} |u|^2 \, dx }   \\
               & {\dm \le \frac4{e(z)^2\mu_0}(1 + s)^{-2(2\delta-1)}      
                        \int_{E_s} (1 + r)^{2\delta-2}|\CD_ru|^2 \, dx } \\
               & \ \ \ \ \ \ \ \ \ \ \ \
                     {\dm + \frac{4{M_0}}{e(z)(2\delta-1)\sqrt{\mu_0}}         
                          (1 + s)^{-(2\delta-1)}
                             \M f \M_{\delta}\M u \M_{-\delta}. }
\endsplit                                                     \T (5.12)
$$
     Note that \f e(z)^{-1} \le 1/\sqrt{2|z|} $. Then, by setting \f s = 0 \A 
     in (5.12) and using Theorem 4.1, it follows that there exists a 
     positive constant \f C_1 = C_1(\delta) \A such that
$$
       \M u \M_{-\delta}^2 \le \frac{C_1}{|z|}\M f \M_{\delta}^2. \T (5.13)
$$
     The estimate (5.7) is obtained from (5.12), (5.13) and Theorem 4.1.        
     \ \ $\Vert$                                 
\BP

        The following corollary is obtained easily when Theorems 5.1 is
     combined with Theorem 5.2.
\BP

        {\bf Corollary 5.3.} {\it Suppose that \r{\em Assumption 2.1}$ holds 
     with \f N \ge 3 $. Let \f 1/2 < \delta \le 1 $. Let \f u = R(z)f $, 
     where \f f \in L_{2,\delta}(\RN) \A and \f z = \lambda + i\eta \A with
     \f \lambda \ge 0 $, \f \eta \ne 0 $. Then there exists a positive 
     constant \f C = C(\delta) \A such that 
$$
\split
    {\dm \int_{E_s} } & {\dm (1 + r)^{-2\delta} 
                    \big( |\gr u|^2 + |k|^2|u|^2 \big) \, dx }\\
          & {\dm \le C(1 + \frac1{\sqrt{|z|}} + \frac1{|z|})
                  (1 + s)^{-(2\delta-1)} \M f \M_{\delta}^2 
                                 \q (s \ge 1), } 
\endsplit                                                      \T (5.14)
$$
     where \f E_s = \{ x \in \RN \ : \ |x| > s \} $, and \f \M \ \M_{t} \A is the 
     norm of \f L_{2,t}(\RN) $. The constant \f C = C(\delta) \A is 
     independent of \f f \A and \f z \A satisfying the above conditions 
     and \f s \ge 1 $.}                                                 
\BP

        Now we are in a position to show some estimate of the operator
     norm of \f R(z) $. For \f 0 < c < d < \infty \A a subset
     \f J_{\pm}(c, \, d) \A of \f \bC \A are defined by
$$
\left\{ \split
      {\dm J_+(c, \, d) = \{ \ z = \lambda + i\eta \ : \ c \le \lambda \le d,
                                     \ 0 < \eta \le 1 \ \}, } \\  *[4pt]
      {\dm J_-(c, \, d) = \{ \ z = \lambda + i\eta \ : \ c \le \lambda \le d,
                                     \ -1 \le \eta < 0 \ \}.  }                              
\endsplit \right.                                                  \T (5.15)
$$
     Let \f t \in \Ro $. The weighted Sobolev spaces \f H_t^j(\RN) $, 
     \f j = 1, 2 $, are defined as the completion of \f \Con(\RN) \A by the 
     norms
$$
        \M u \M_{1,t} 
          = \bigg[ \int_{\RN} (1 + r)^{2t} \big( |\gr u|^2 + 
                                |u(x)|^2 \big) \, dx \bigg]^{1/2},  \T (5.16)
$$
     and
$$
        \M u \M_{2,t} 
          = \bigg[ \int_{\RN} (1 + r)^{2t} \sum_{|\gamma|\le 2}
                        |\pa^{\gamma}u|^2 \, dx \bigg]^{1/2},      \T (5.17)                           
$$
     respectively, where
$$
\left\{ \split
       {\dm \gamma = (\gamma_1, \gamma_2, \cdots, \gamma_N), } \\ *[4pt]
       {\dm |\gamma| = \gamma_1 + \gamma_2 + \cdots + \gamma_N, } \\ *[4pt]
       {\dm \pa^{\gamma}u = (\pa_1)^{\gamma_1} \cdots (\pa_N)^{\gamma_N}u
            \q (\pa_j = \pa/\pa x_j). }
\endsplit \right.                                                 \T (5.18)
$$
    The inner product and norm of \f H_t^j(\RN) \A will be denoted by
    \f (\ , \ )_{j,t} \A and \f \M \ \M_{j,t} $. For an operator \f T $,\, 
    the operator norm in \f \bB(H_s^j(\RN), \, H_t^{\ell}(\RN)) $ will
    be denoted by  \f \M T \M_{(j,s)}^{(\ell,t)} $, where \f j, \ell = 0, 1, 
    2 $, \f s, t \in \Ro $, and we set
$$
         H_s^0(\RN) = L_{2,s}(\RN).                                 \T (5.19)
$$         
\BP

        {\bf Theorem 5.4.} {\it Suppose that \r{\em Assumption 2.1}$ holds
     with \f N \ge 3 $. Let \f 1/2 < \delta \le 1 \A. Let \f R(z) \A be the 
     resolvent of \f H $.        
\SP

        \r{\em (i)}$ Then there exists a positive constant 
     \f C = C(\delta) \A such that 
$$
         \M R(z) \M_{(0,\delta)}^{(0,-\delta)} \le \frac{C}{\sqrt{|z|}},    
                                                                  \T (5.20)    
$$
     for \f z = \lambda + i\eta \in \bC \A with \f \lambda \ge 0 \A and
     \f \eta \ne 0 $.
\SP

        \r{\em (ii)}$ Let \f 0 < c < d < \infty \A and let 
     \f J_{\pm}(c, \, d) \A be as above. Then there exists a positive 
     constant \f C = C(\delta, \, c, \, d) \A such that
$$
         \M R(z) \M_{(0,\delta)}^{(2,-\delta)} \le C              \T (5.21)    
$$             
     for \f z \in J_+(c, \, d)\cup J_-(c, \, d) $.}
\BP

        {\it Proof. \ } (i) directly follows from (5.10) in Theorem 5.3
     with \f s = 0 $. It follows from (4.2) in Theorem 4.1 that
$$
\split
      {\dm \int_{E_1} } & {\dm (1 + r)^{-2\delta} |\gr u|^2 \, dx } \\ *[6pt]
           & {\dm \le 2 \int_{E_1} (1 + r)^{-2\delta} |\CD u|^2 \, dx } 
                                                                    \\ *[6pt]
           &  \hspace{3cm} {\dm  + 2 \int_{E_1} (1 + r)^{-2\delta} 
                    \bigg|\frac{N-1}{2r}\tx u - ik\tx u \bigg|^2 \, dx }  
                                                                    \\ *[6pt]
           & {\dm \le 2 \int_{E_1} (1 + r)^{2\delta-2} |\CD u|^2 \, dx }  
                                                                    \\ *[6pt]                                       
           &  \hspace{3cm} 
                 {\dm + 4 \int_{E_1} \bigg( \big(\frac{N-1}{2}\big)^2
                   + |k|^2 \bigg) (1 + r)^{-2\delta} |u|^2 \, dx } \\  *[6pt]
                 & {\dm \le C_2(\delta) \M f \M_{\delta}^2 +
                            C_3(\delta, \, c, \, d) \M u \M_{-\delta}^2 }
\endsplit                                                       \T (5.22)
$$
     with positive constants \f C_2 = C_2(\delta) \A and 
     \f C_3 = C_3(\delta, \, c, \, d) $, where \f u = R(z)f \A \linebreak 
     with \f f \in L_{2,\delta}(\RN) \A and
     \f z \in J_+(c, \, d)\cup J_-(c, \, d) $. Since the \f L_2(B_1) $-norm
     of \f \gr u \A can be evaluated by the interior estimate, we 
     obtain from (5.22) and (5.20)                                  
$$
     \M u \M_{1,-\delta} \le C_4(\delta, \, c, \, d)\M f \M_{2,\delta}
                                                               \T (5.23) 
$$
     with a positive constant \f C_4 = C_4(\delta, \, c, \, d) $, which 
     implies that
$$
    \M R(z) \M_{(0,\delta)}^{(1,-\delta)} \le C_5(\delta, \, c, \, d) 
                         \qq          (z \in J_+(c, \, d)\cup J_-(c, \, d))
                                                             \T (5.24)
$$
     with a constant \f C_5 = C_5(\delta, \, c, \, d) $. Using the relation
     \f \Delta u = - \,\mu f - k^2u $, we have from (5.20)                                              
$$
     \M \Delta u \M_{-\delta} \le C_6(\delta, \, c, \, d)\M f \M_{\delta}
                                                                  \T (5.25)
$$
     with a constant \f C_6 = C_6(\delta, \, c, \, d) $. The inequality 
     (5.21) follows from (5.23) and (5.24) ([\JS], Proposition A.2), which 
     completes the proof. \ \ $\Vert$

$$
\ \ \ \ \
$$

\cen{     {\bf \S6. Limiting absorption principle.}}           
\SP

        By the use of the results established in \S3, \S4 
     and \S5, we can show the limiting absorption principle for the 
     operator \f H \A in \f \RN \A with \f N \ge 3 \A using the 
     arguments used to prove the limiting 
     absorption principle for the \Sch\ operator (e.g., [\Sa], [\IS]).
     
        First we shall define the boundary value \f R^{\pm}(\lambda) $, 
     \f \lambda > 0 $, of the resolvent \f R(z) \A when 
     \f z = \lambda + i\eta \to \lambda $.
\BP

        {\bf Theorem 6.1.} {\it Suppose that \r{\em Assumption 2.1}$ holds 
     with \f N \ge 3 $. Let \f 1/2 < \delta \le 1 \A. Let 
     \f f \in L_{2,\delta}(\RN) \A and let \f \lambda > 0 $. Then there 
     exist
$$
\left\{ \split
        {\dm  \lim_{\eta\downarrow 0} R(\lambda + i\eta)f 
                              = u_{+}(\cdot, \, \lambda, \, f), } \\ *[4pt]
        {\dm  \lim_{\eta\downarrow 0} R(\lambda - i\eta)f 
                              = u_{-}(\cdot, \, \lambda, \, f), }\\
\endsplit \right.                                             \T (6.1)
$$
     in \f H_{-\delta}^2(\RN) $, where \f H_{-\delta}^2(\RN) \A is given in
     \S 5, and \f u_{+}(\cdot, \lambda, f) \A \, \r{\em [}$\,or  
     \f u_{-}(\cdot, \lambda, f) $\,\r{\em ]}$ is a unique solution of the
     equation 
$$
\left\{ \split
       {\dm  - \mu(x)^{-1} \Delta u - \lambda u = f, } \\ *[4pt]
       {\dm \M \CD^{(+)}u \M_{\delta-1} < \infty \ \ \ \
            [ \mbox{\rm or \ } \M \CD^{(-)}u \M_{\delta-1} < \infty \,] }
\endsplit \right.                                                 \T (6.2)
$$        
     with}
$$
       \CD^{(\pm)}u = \gr u + \{(N-1)/(2r)\}\tx u 
                               \mp i\sqrt{\lambda\mu(x)}\tx u.     \T (6.3)    
$$
\BP
                      
        {\it Proof. \ } For each \f n = 1, 2, \cdots \A let 
     \f z_n = \lambda + i\eta_n $, where \f \eta_n > 0 \A and 
     \f \eta_n \downarrow 0 \A as \f n \to \infty $. Set \f u_n = R(z_n)f $,
     \f n = 1, 2, \cdots $. Then, in view of Corollary 5.3 and Theorems 5.4, 
     we see that not only the sequence \f \{ u_n \} \A is a bounded set in
     \f H_{-\delta}^2(\RN) \A but also \f \M u_n \M_{1, -\delta, E_s} \A 
     is uniformly small for \f n \A as \f s \to \infty $, where
$$
      \M u_n \M_{1, -\delta, E_s}^2 
             = \int_{E_s} (1 + |x|)^{-2\delta}|u_n|^2 \, dx.        \T (6.4)
$$
     Therefore, by the Rellich selection theorem, \f \{ u_n \} \A has a
     subsequence which converges to a limit function \f u_0 \A in     
     \f H_{-\delta}^2(\RN) \A (see, e.g., [\JS], Proposition A.3). Since
     \f u_0 \A turns out to be a unique solution of the equation
     \f - \, \Delta u - k^2u = \mu f \A with radiation condition
     \f \M \CD^{(+)}u \M_{\delta-1} < \infty $, i.e., 
     \f u_0 =  u_{+}(\cdot, \, \lambda, \, f) $, it follows that the 
     sequence \f \{ u_n \} \A itself converges to 
     \f u_0 =  u_{+}(\cdot, \, \lambda, \, f) \A in \f H_{-\delta}^2(\RN) $.
     The existence of the first limit of (6.1) follows from the above 
     argument. The existence of the second limit of (6.1) 
     can be proved in the same way. \ \ $\Vert$
\BP

        {\bf Definition 6.2.} Let \f \lambda > 0 $. Then the operators
     \f R_{\pm}(\lambda) \A are defined by
$$
\left\{ \split
      R_+(\lambda) \, : \, L_{2,\delta}(\RN) \ni f \longmapsto 
                             u_{+}(\cdot, \, \lambda, \, f) \in              
                                             H_{-\delta}^2(\RN), \\ *[4pt]
      R_-(\lambda) \, : \, L_{2,\delta}(\RN) \ni f \longmapsto 
                             u_{-}(\cdot, \, \lambda, \, f) \in              
                                             H_{-\delta}^2(\RN), \\
\endsplit \right.                                              \T (6.5)
$$
\BP
                                              
        Let \f D_{\pm} \subset \bC \A be given by
$$
\left\{ \split
  D_+ 
     = \{ \, z = \lambda + i\eta \, : \, \lambda > 0, \, \eta \ge 0 \, \}, 
                                                                  \\ *[4pt]       
  D_-
     = \{ \, z = \lambda + i\eta \, : \, \lambda > 0, \, \eta \le 0 \, \}.     
\endsplit \right.                                               \T (6.6)
$$
     Then the resolvent \f R(z) \A will be extended on each of \f D_{\pm} \A
     by the use of \f R_{\pm}(\lambda) $, i.e., for \f z \in D_+ \A we set
$$
     R(\lambda+i\eta) = 
\left\{ \split
    R(\lambda+i\eta) \qq (\lambda > 0, \, \eta > 0), \\  *[4pt]      
    R_+(\lambda) \ \ \qq \ \ (\lambda > 0, \, \eta = 0),     
\endsplit \right.                                               \T (6.7)
$$
     and for \f z \in D_- \A we set                  
$$
     R(\lambda+i\eta) = 
\left\{ \split
    R(\lambda+i\eta) \qq (\lambda > 0, \, \eta < 0), \\ *[4pt]       
    R_-(\lambda) \ \ \qq \ \ (\lambda > 0, \, \eta = 0).     
\endsplit \right.                                               \T (6.8)
$$
     For \f 0 < c < d < \infty \A let \f J_{\pm}(c, d) \A
     be as in (5.15). The closure  \f \overline{J}_{\pm}(c, d) \A 
     \linebreak are given by
$$
\left\{ \split
      \overline{J}_+(c, \, d) = \{ \ z = \lambda + i\eta \ : \ 
          c \le \lambda \le d, \ 0 \le \eta \le 1 \ \} \subset D_+, \\ *[4pt]   
      \overline{J}_-(c, \, d) = \{ \ z = \lambda + i\eta \ : \ 
          c \le \lambda \le d, \ - \, 1 \le \eta \le 0 \ \} \subset D_-.                                
\endsplit \right.                                                  \T (6.9)
$$     
        For \f \lambda \in D_+ \cap (0, \, \infty) \A [or 
     \f D_- \cap (0, \, \infty)$], \f \CD u \A should be interpreted 
     as \f \CD^{(+)} \A [\,or \f \CD^{(-)} $].    
          
        From Theorems 6.1, 5.1, 5.2, 5.4 and Corollary 5.3, we easily see
     the following:   
\BP

        {\bf Theorem 6.3.} {\it Let \f 1/2 < \delta \le 1 $. Suppose
     \r{\em Assumption 2.1}$ holds with \f N \ge 3 $. Let \f R(z) \A be 
     extended on each of \f D_+ \A and \f D_- $. 
\SP
             
        \r{\em (i)}$ Then there exists a positive constant 
     \f C = C(\delta) \A such that 
$$
\left\{ \split
     {\dm \int_{E_s} (1 + r)^{-2\delta} |R(z)f|^2 \, dx 
         \le \frac{C^2}{|z|}(1 + s)^{-2(2\delta-1)} \M f \M_{\delta}^2 } 
                                                                   \\ *[5pt]
     {\dm \ \ \ \ \ \ \ \ \ \ \ \ \ \ \ \ \ \ \ \ \ \ \ \ \ \ \ \ \ \ \ \
           \ \ \ \ \ \ \ \ \ \ \ \ \ \  
                            (s \ge 0, \ f \in L_{2,\delta}(\RN)),} \\ *[5pt] 
     {\dm \M R(z) \M_{(0,\delta)}^{(0,-\delta)} \le \frac{C}{\sqrt{|z|}}
                                    \q (z \in D_+ \cup D_-), } \\ *[5pt]
     {\dm \M \CD R(z)f \M_{\delta-1} \le C \M f \M_{\delta} \q
             (z \in D_+ \cup D_-, \ f \in L_{2,\delta}(\RN)). }                            
\endsplit \right.                                              \T (6.10)
$$      
\SP

        \r{\em (ii)}$ For \f 0 < c < d < \infty \A there exists a positive 
     constant \f C = C(c, \ d, \ \delta) \A \linebreak such that, for 
     \f z \in \overline{J}_+(c, \, d) \cup \overline{J}_-(c, \, d) $,}
$$   
\left\{ \split
      {\dm \M R(z) \M_{(0,\delta)}^{(2,-\delta)} \le C, } \\ *[5pt]
      {\dm \int_{E_s}(1 + r)^{-2\delta} 
                \big( |\gr R(z)f|^2 + |k|^2|R(z)f|^2 \big) \, dx } \\  *[5pt]
      \hspace{3.5cm}   \le C^2 (1 + s)^{-(2\delta-1)} \M f \M_{\delta}^2   
                                                                  \\  *[4pt]
       \ \ \ \ \ \ \ \ \ \ \ \ \ \ \ \ \ \ \ \ \ \ \ \ \ \ \ \ \ \ \ \ \
         \ \ \ \ \ \ \ \  (s \ge 1, \ f \in L_{2,\delta}(\RN)).  
\endsplit \right.                                              \T (6.11)
$$
\BP

        The next proposition will be used when we prove continuity of
     \f R(z) \A with respect to \f z \A and the compactness of the operator
     \f R(z) $.   
\BP

        {\bf Proposition 6.4.} {\it Let \f 1/2 < \delta \le 1 $. Suppose
     that \r{\em Assumption 2.1}$ holds with \f N \ge 3 $. Let \f R(z) \A be 
     extended on each of \f D_+ \A and \f D_- $. Let \f \{ f_n \} \A be a 
     sequence in \f L_{2,\delta}(\RN) \A such that
$$
             f_n \to f_0 \q {\rm weakly \ \, in \ \, } 
                                               L_{2,\delta}(\RN)    \T (6.12)
$$
     as \f n \to \infty $, and let \f \{ z_n \} \subset 
     \overline{J}_+(c, \, d) \A {\rm [}\,or \f \{ z_n \} \subset 
     \overline{J}_+(c, \, d) $\,{\rm ]} with \f 0 < c < d < \infty \A such 
     that 
$$
            z_n \to z_0 \qq (n \to \infty).                   \T (6.13)
$$
     Then there exists a sequence \f \{ n_k \}_{k=1}^{\infty} \A of positive
     integers such that 
$$
        n_1 < n_2 < n_3 < \cdots < n_k < \cdots \to \infty,   \T (6.14)
$$
     and 
$$
             u_{n_k} \to R(z_0)f_0 \q {\rm \ in \ } H_{-\delta}^1(\RN)
                                                              \T (6.15)
$$ 
     as \f k \to \infty $, where \f u_{n_k} = R(z_{n_k})f_{n_k} $.}                                                                       
\BP

        {\it Proof. \ } We are going to give the proof for the case that
     \f z_0 = \lambda_0 \in [c, \, d] \A and \f \{ z_n \} \subset 
     \overline{J}_+(c, \, d) $. The case that \f z_0 \A is not a 
     real number can be treated more easily. Set \f u_n =  R(z_n)f_n $.
     Using Theorem 6.3 and proceeding as in the proof
     of Theorem 6.1, we see that there exists a subsequence of 
     \f \{ u_n \} \A which converges to a limit function \f u_0 \A in 
     \f H_{-\delta}^1(\RN) $. Then it is easy to show that 
     \f u_0 \in H_{-\delta}^2(\RN) \A and that
     \f u_0 = u_+(\cdot, \lambda_0, f_0) = R(\lambda_0)f_0 $, which 
     completes the proof. \ \ $\Vert$    
\BP

         The following properties of the extended resolvent \f R(z) \A 
     follows directly from the above proposition.
\BP

        {\bf Theorem 6.5.} {\it Let \f 1/2 < \delta \le 1 $. Suppose
     that \r{\em Assumption 2.1}$ holds with \f N \ge 3$. Let \f R(z) \A be 
     extended on each of \f D_+ \A and \f D_- $.    
\SP

        \r{\em (i)}$ Then the extended resolvent \f R(z) \A is a 
     \f \bB(L_{2,\delta}(\RN), \, H_{-\delta}^2(\RN))$-valued continuous 
     function on each of \f D_+ \A and \f D_- $. 
\SP

        \r{\em (ii)}$ For any \f z \in D_+ \A $\mb{\em [}$or $ D_- \,
     \mb{\em ]}$, \f R(z) \A is a compact operator from \linebreak 
     \f L_{2,\delta}(\RN) \A into \f H_{-\delta}^1(\RN) $.}      
\BP

        {\it Proof. \ } (I) The proof of (i). Suppose that there is 
     \f z_0 \in D_+ \A at which \f R(z) \A is not continuous in the topology 
     of \f \bB(L_{2,\delta}(\RN), \, H_{-\delta}^1(\RN)) $. We may assume 
     that \f z_0 = \lambda_0 > 0 $, since the other case can be handled 
     more easily. Then there exist \f \epsilon_0 > 0 \A and sequences 
     \f \{ z_n \} \subset D_+ $, \f \{ f_n \} \subset L_{2,\delta}(\RN) \A 
     and \f f_0 \in L_{2,\delta}(\RN) \A such that
$$     
\left\{ \split
      z_n \to \lambda_0 \qq \hspace{1.3cm} (n \to \infty), \\ *[4pt]
      \M f_n \M_{\delta} = 1 \qq \hspace{1cm} (n = 1, 2, \cdots), \\  *[4pt]
      f_n \to f_0 \q {\rm weakly \ \, in \ \, } 
                                               L_{2,\delta}(\RN), \\  *[4pt]    
         \M R(\lambda_0)f_n - R(z_n)f_n \M_{1,-\delta} \ge \epsilon_0.
\endsplit \right.                                             \T (6.16)
$$
     Applying Proposition 6.4 for the sequence
     \f \{ R(\lambda_0)f_n \} $, we see that there exists a subsequence
     \f \{ f_{n_k} \} \A of \f \{ f_n \} \A such that 
$$
            R(\lambda_0)f_{n_k} \to R(\lambda_0)f_0  \qq      
                            {\rm in \ } \  H_{-\delta}^1(\RN)      \T (6.17)
$$
     as \f k \to \infty $. Apply Proposition 6.4 again for the sequence 
     \f \{ R(z_{n_k})f_{n_k} \}_{k=1}^{\infty} \A to see that there is a 
     subsequence \f \{ R(z_{n_{k_p}})f_{n_{k_p}} \}_{p=1}^{\infty} \A
     such that
$$
      R(z_{n_{k_p}})f_{n_{k_p}} \to R(\lambda_0)f_0 \qq  
                            {\rm in \ } \  H_{-\delta}^1(\RN).      \T (6.18)
$$
     as \f p \to \infty $. Therefore it follows from (6.17) and (6.18) that
$$
\split
   & \M R(\lambda_0)f_{n_{k_p}} - R(z_{n_{k_p}})f_{n_{k_p}} \M_{1,-\delta} 
                                                                 \\  *[4pt]
   & \hspace{2cm} \le \M R(\lambda_0)f_{n_{k_p}} 
                              - R(\lambda_0)f_0 \M_{1,-\delta} \\  *[4pt]
   & \hspace{2cm} \ \ \ \ \ \ \ \ \ \ \ \ \ \ \ \ \
   + \M R(\lambda_0)f_{0} - R(z_{n_{k_p}})f_{n_{k_p}} \M_{1,-\delta} 
                                                                 \\ *[4pt] 
   & \hspace{2cm} \to 0 
\endsplit                                                          \T (6.19) 
$$   
     as \f p \to \infty $, which contradicts the fourth relation of (6.16). 
     Thus we have shown that \f R(z) \A is a 
     \f \bB(L_{2,\delta}(\RN), \, H_{-\delta}^1(\RN))$-valued 
     continuous function on each of \f D_+ \A and \f D_- $. 
\SP
 
        (II) Proof of (i) (continued). Let \f z, \, z_0 \in D_+ \A [or
     \f z, \, z_0 \in D_- $]. The continuity of \f R(z) \A in
     \f \bB(L_{2,\delta}(\RN), \, H_{-\delta}^1(\RN)) \A is combined with
     the relation
$$
\split
    & \Delta R(z) - \Delta R(z_0) \\  *[4pt]
    & \hspace{2cm}           = - z\mu R(z) + z_0\mu R(z_0)  \\ *[4pt]
    & \hspace{2cm}           = (z_0 - z)R(z_0) 
                         +  z(R(z_0) - R(z)) \to 0      
\endsplit                                                          \T (6.20)
$$
     in \f \bB(L_{2,\delta}(\RN), \, L_{2,-\delta}(\RN)) \A as 
     \f z \to z_0 $ to obtain the continuity of \f R(z) \A in
     \f \bB(L_{2,\delta}(\RN), \, H_{-\delta}^1(\RN)) \A (cf., e.g.,
     [\JS], Proposition A.3 in Appendix A.2). This completes the proof of (i).
\SP

        (III) Proof of (ii). Let \f \{ f_n \} \A be a bounded sequence
     in \f L_{2,\delta}(\RN) $. We may assume with no loss of generality
     that the sequence \f \{ f_n \} \A converges weakly in      
     \f L_{2,\delta}(\RN) $. The weak limit will be denoted by \f f_0 $.
     Then, applying Proposition 6.4, we see that there exists a subsequence
     \f \{ f_{n_k} \}_{k=1}^{\infty} \A such that
$$
          R(z)f_{n_k} \to R(z)f_0 \q {\rm \ in \ } H_{-\delta}^1(\RN)
                                                              \T (6.21)
$$ 
     as \f k \to \infty $, which completes the proof of (ii). \ \ $\Vert$
\BP

        (i) of Theorem 6.5 and the spectral formula for self-adjoint 
     operators are combined to give
\BP

        {\bf Corollary 6.6} {\it Suppose that \r{\em Assumption 2.1}$ holds 
    with \f N \ge 3 $. Then the selfadjoint operator \f H \A is absolutely 
    continuous on the interval \f (0, \infty) $.}
                   
$$
\ \ \ \ \ 
$$

\cen{  {\bf \S7. The operator \f H \A in \f \Rb $.} }              
\BP

        In the two dimensional case, the constant \f c_N \A given by (3.8)
     takes the value \f c_2 = -1/4 < 0 \A although \f c_N \ge 0 \A for all
     \f N \ge 3$. Because of this, we are going to make some technical 
     changes in the theory which was dveloped in \S3 $\sim$ \S6. Also we 
     should note that \f u/|x| \A is not necessarily integrable around 
     \f x = 0 \A for \f u \in H^2(\Rb)_{{\rm loc}} \A although we have 
$$
        \frac{u}{|x|^{1/2}} \in L_2(\Rb)_{{\rm loc}}   
                              \q (u \in H^2(\Rb)_{{\rm loc}})    \T (7.1)
$$
     since  \f u \in H^2(\Rb)_{{\rm loc}} \A is a continuous function on
     \f \Rb $. We are going to use Notation 3.1 (with \f N = 2 $) throughout
     this section.
\BP

        {\bf 7.1.} \underline{Uniqueness of the solution}.
\SP

       The uniqueness theorem takes the following form:
\BP

       {\bf Theorem 7.1.}\ {\it Suppose that \r{\em Assumption 2.1}$ with 
     \f N = 2 \A holds. Let \linebreak \f u \in H^2(\RN)_{{\rm loc}} \A be a 
     solution of the homogeneous equation 
$$
           -\, \mu(x)^{-1}\Delta u - \lambda u = 0 \qq (\lambda > 0)
                         					                         \T (7.2)
$$
     on \f \Rb \A such that                           					
$$
     \liminf_{R\to\infty} 
         R^{\alpha}\int_{S_R} \big( \bigg|\frac{\pa u}{\pa r}\bigg|^2 
                        + |u|^2 \big)\, dS = 0                    \T (7.3)
$$                        
     with \f \alpha > 0$. Then \f u \A is identically zero.}
\BP

     {\it Proof. \ } (I) Note that Proposition 3.3 and Lemma 3.4 are true
     for the case of \f N = 2 \A without any change. Also Lemma 3.5 is true
     if we add a condition that
$$
             \frac{\pa \varphi}{\pa r} = O(r)                      \T (7.4)
$$
     as \f r \downarrow 0$.

        (II) We may assume with no loss of generality that 
     \f 0 < \alpha \le 1 $. Let
$$
       \varphi(x) =
\left\{ \split
                 |x|^2  \qq \q (0 \le |x| \le r_0), \\  *[4pt]
                 r_0^{2-\alpha}|x|^{\alpha} \qq (|x| > r_0)
\endsplit \right.                                                  \T (7.5)          
$$    
     in (3.7), where \f r_0 > 0 \A will be determined later. Proceeding
     as in the proof of Theorem 3.2, we have from (3.7), for any 
     \f R > r_0 > r > 0 $, 

\newpage

$$
\split
  &{\dm  \int_{B_{R}} \frac12 \frac{\pa \varphi}{\pa r} k^2|u|^2\, dx
     + \int_{B_{rR}} \big(\frac{\varphi}r - 2^{-1}\frac{\pa\varphi}{\pa r}
        \big (|\gr u|^2 - \bigg|\frac{\pa u}{\pa r}\bigg|^2)\, dx } \\ *[6pt]
  &  \ \ \ \ \hspace{5cm}
           {\dm   - \frac14 \int_{B_{rR}} r^{-2}\big( \frac{\varphi}r 
              - 2^{-1}\frac{\pa\varphi}{\pa r} \big) |u|^2\, dx  } \\  *[6pt]
  &  \hspace{0.5cm} {\dm  \le 2^{-1}\sum_{\ell=1}^2 
                         \int_{\pa\Omega_{\ell}\cap B_{rR}}
                       \varphi |k|^2(\tx \cdot n)|u|^2\, dS } \\ *[6pt]
  &   \ \ \ \  \hspace{0.5cm} {\dm + 2^{-1} r_0^{2-\alpha} R^{\alpha}
                \int_{S_R} \bigg(2|\CD_ru|^2 - |\CD u|^2 + 
                   \frac14 r^{-2}|u|^2 
              {\rm Im}\big( k\frac{\pa u}{\pa n} \overline{u} \big) \bigg) \, dS }  
                                                                  \\  *[6pt]              
  &  \ \ \ \ \hspace{0.5cm} {\dm - 2^{-1} r^2\int_{S_r} (2|\CD_ru|^2 
            - |\CD u|^2 + \frac14 r^{-2}|u|^2)\, dS, }  \\  *[6pt]
  &  \ \ \ \ \hspace{0,5cm} {\dm  + \int_{B_r} \frac12 r^2 |\CD u|^2\, dx 
          + \sum_{\ell=1}^{2} \int_{\pa\Omega_{\ell}\cap B_r}
                r^2 {\rm Im}\big( k\frac{\pa u}{\pa r}
                                 \overline{u}\big) \, dS.  }    
\endsplit                                                          \T (7.6)                                               
$$
     The left-hand side and right-hand side of (7.6) will be denoted by 
     \f K_L \A and \f K_R $, respectively.

        (III) As in the proof of Theorem 3.2, we see that the second term of 
     \f K_L \A is nonnegative: Thus we have   
$$
\split
    & {\dm K_L \ge \int_{B_{r_0}} r k^2|u|^2\, dx } \\ *[5pt]
    & \hspace{1.5cm} {\dm + \int_{B_{r_0R}} 
        \bigg( \frac12 \frac{\pa \varphi}{\pa r} k^2
      - \frac14 r^{-2} \big( \frac{\varphi}r - 2^{-1}\frac{\pa\varphi}{\pa r} \big)
                                     \bigg) |u|^2\, dx } \\ *[5pt]
    & \hspace{0.5cm} {\dm = \int_{B_{r_0}} r k^2|u|^2\, dx } \\ *[5pt]
    & \hspace{2cm} {\dm + r_0^{2-\alpha}
                 \int_{B_{r_0R}} r^{\alpha-1} \bigg( \frac{\alpha}2  k^2
           - \frac14 \big( 1 - \frac{\alpha}2 \big) r^{-2} \bigg) |u|^2\, dx. } 
\endsplit                                                           \T (7.7)  
$$     
     Choose \f r_0 = r_0(k, \alpha) \A so large that
$$
      \frac{\alpha}2  k^2 - \frac14 \big( 1 - \frac{\alpha}2 \big) r^{-2}   
                        \ge \frac{\alpha k^2}4  \q ( r \ge r_0 )    \T (7.8)
$$
     Then we have
$$
       K_L \ge  \int_{B_{r_0}} r k^2|u|^2\, dx + 
           r_0^{2-\alpha} \int_{B_{r_0R}} r^{\alpha-1} 
                       \frac{\alpha k^2}4 |u|^2\, dx.             \T (7.9)
$$

        (IV) We cn proceed as in the proof Theorem 3.2 to see that the first 
     term of \f K_R \A is nonpositive, the fourth and fifth term go to 
     \f 0 \A as \f r \downarrow 0 $, and the third term goes to \f 0 \A as 
     \f r \downarrow 0 \A along an appropriate sequence. Therefore, after 
     evaluating the third term using \f |u| \A and \f |\pa u/\pa u| $, we 
     obtain 
$$
\split
   & {\dm  \int_{B_{r_0}} r k^2|u|^2\, dx + 
    r_0^{2-\alpha} \int_{B_{r_0R}} r^{\alpha-1} \frac{\alpha k^2}4 |u|^2\, dx } 
                                                                \\ *[7pt]
   & \hspace{2cm} 
      {\dm \le CR^{\alpha} \int_{S_R} \big( \bigg|\frac{\pa u}{\pa r}\bigg|^2 
                                    + |u|^2 \big)\, dS }                                                                         
\endsplit                                                           \T (7.10)
$$
     with \f C = C(\lambda, \alpha) $, and hence it follows from the 
     condition (7.3) that 
$$
       \int_{B_{r_0}} r k^2|u|^2\, dx + 
      r_0^{2-\alpha} 
         \int_{E_{r_0}} r^{\alpha-1} \frac{\alpha k^2}4 |u|^2\, dx = 0,
                                                                   \T (7.11)
$$
     where \f E_{r_0} = \{ \ x \in \Rb \ : \ |x| > r_0 \ \} $, i.e., \f u \A 
     is identically zero, which completes the proof. \ \ $\Vert$
\BP

        Since we have established Theorem 7.1, a two dimensional version of 
     Theorem 3.2, we can easily see that each of Corollary 3.6, Theorem 3.7 
     and Corollary 3.8 has its two dimensional version only by replacing the 
     conditions (3.24), (3.25), and (3.35) by
$$
      \liminf_{R\to\infty} R^{\alpha}\int_{S_R} |\CD_r^{(\pm)} u|^2\, dS = 0,    
                                                                   \T (7.12)
$$
$$
     \liminf_{R\to\infty} R^{\alpha}\int_{S_R} |\frac{\pa u}{\pa r} 
                                        \mp iku|^2\, dS = 0,        \T (7.13)
$$       
     and
$$
\left\{ \split                                                       
   {\dm \int_{E_R} r^{-1+\alpha} |\CD_r^{(\pm)} u|^2 \, dx < \infty, }  \\  
                                                                     *[4pt]
   {\dm \int_{E_R} r^{-1+\alpha} |\frac{\pa u}{\pa r} \mp iku|^2, dx 
                                                              < \infty, }
\endsplit \right.                                                 \T (7.14)       
$$      
     with \f \alpha > 0 $. We do not take the trouble to write down these
     two dimensional versions since they are now quite obvious.
\BP

        {\bf 7.2.} \underline{The evaluation of \f \CD u $}.
\SP

        Consider \f u \A given by
$$
\left\{ \split
          u = R(z)f,  \\ *[4pt]
          z = \lambda + i\eta  \qq (\lambda \ge 0, \eta \ne 0),  \\ *[4pt]
          f \in L_{2, \delta}(\Rb).
\endsplit \right.                                            \T (7.15)
$$
\BP

        {\bf Theorem 7.2.} {\it Assume \r{\em Assumption 2.1}$ with 
     \f N = 2$. Let \f 1/2 < \delta \le 1 $. Let \f 0 < c < d < \infty \A 
     and let \f J_{\pm}(c, d) \A be as in \r{\em (5.15)}$. Let \f u \A be 
     given by \r{\em (7.15)}$ with \f z \in J_{+}(c, d) \cup J_{-}(c, d) $. 
     Then there exists a positive constant \f C = C(c, d, \mu, \delta) \A 
     such that
$$
        \M \CD u \M_{\delta-1. *} 
                  \le C \big( \M f \M_{\delta} + \M u \M_{-\delta} \big)
                                                                 \T (7.16)
$$
     where }
$$
       \M v \M_{t, *}^2 = \int_{B_1} |x||v(x)|^2 \, dx
                    + \int_{E_1} (1 + |x|)^{2t}|v(x)|^2 \, dx.   \T (7.17)
$$
\BP

        {\it Proof. \ } Set \f \varphi(x) = \xi(|x|)/\sqrt{\mu(x)} \A in 
     (3.7), where
$$
      \xi(r) =
\left\{ \split
           {\dm \frac12 r^2 \qq \hspace{1.9cm} (r \le 1/2),} \\ *[4pt]
           {\dm \frac{1}{2^{2\delta}}(1 + r)^{2\delta-1}   \qq (r \ge 1) }.
\endsplit \right.                                                  \T (7.18)
$$
     We can evaluate each term in (3.7) in quite a similar manner as in the
     Proof of theorem 4.1 except the fourth term \f I_{L4} \A of the 
     left-hand side which is nonpositive in our case because 
     \f c_2 = -1/4 < 0 $. The term \f - I_{L4} \A can be evaluated as
$$
\split
     & {\dm - I_{L4} = \frac14 \int_{B_{1/2, R}} r^{-2} \bigg(
                     \frac{\varphi}{r} - 2^{-1}\frac{\pa \varphi}{\pa r}
                            + b\varphi \bigg) |u|^2 \, dx } \\ *[6pt] 
     & \hspace{1cm} {\dm \le C_1 \M u \M_{-\delta}^2 
                          + C_2 \int_{\Rb} |\eta||u|^2 \, dx, } 
\endsplit                                                       \T (7.19)
$$
     and the second term of the right-hand side of (7.19) is eavaluated as
$$
    \int_{\Rb} |\eta||u|^2 \, dx \le C_3(|f|, |u|)_0            \T (7.20)
$$
     (see, e.g., Eidus\,[\Ei], [\Sd], Lemma 2.1), where 
     \f C_1 = C_1(\mu_0, \delta) $, \f C_2 = C_2(c, d) $, and
     \f C_3 = C_3(c, d, \mu) $. Thus, using (7.19) and (7.21), we can proceed as 
     in the proof of Theorem 4.1 to obtain (7.16), 
     which completes the proof. \ \ $\Vert$
\BP

        The following corollary is now obvious.
\BP

        {\bf Corollary 7.3.} {\it Let \f u = R(z)f \A be ae in 
     {\rm Theorem 7.2.} Then there exists a positive constant 
     \f C = C(c, d, \mu, \delta) \A such that} 
$$
      \int_{E_1} (1 + |x|)^{2\delta-2}|\gr u - ik\tx u|^2 \, dx
                \le C \big( \M f \M_{\delta}^2 + \M u \M_{-\delta}^2 \big).
                                                                  \T (7.21) 
$$

\newpage

        {\bf 7.3.} \underline{Boundedness of \f R(z) \A and the limiting 
     absorption principle}.
\SP

        The following theorem can be proved in quite the same manner as 
     in the proof of Theorem 5.1.
\BP

        {\bf Theorem 7.4.} {\it Suppose that \r{\em Assumption 2.1}$ holds
     with \f N = 2 $. Let \f 1/2 < \delta \le 1 $. Let 
     \f 0 < c < d < \infty \A and let \f J_{\pm}(c, d) \A be as in 
     \r{\em (5.15)}$. Let \f u \A be given by \r{\em (7.15)}$ with 
     \f z \in J_{+}(c, d) \cup J_{-}(c, d) $. Then there exists a positive 
     constant \f C = C(c, d, \mu, \delta) \A such that
$$
    \int_{E_s} (1 + r)^{-2\delta} 
                    \big( |\gr u|^2 + |k|^2|u|^2 \big) \, dx  \hspace{6cm}
$$
$$
          \le C(1 + s)^{-(2\delta-1)} 
                    \big( \M f \M_{\delta}^2 + \M u \M_{-\delta}^2 \big)
                                 \q (s \ge 1),  
                                                      \T (7.22)
$$
     where \f E_s = \{ x \in \Rb \ : \ |x| > s \} $, and \f \M \ \M_{t} \A 
     is the norm of \f L_{2,t}(\Rb) $.}                         
\BP

        Let \f H_t^j(\Rb) $, \f j = 1, 2 $, be defined in \S5 with 
     \f N = 2 $. In order to obtain the counterpart of Theorem 5.4, we 
     prepare 
\BP

        {\bf Proposition 7.5.} {\it Assume \r{\em Assumption 2.1}$ with 
     \f N = 2$. Let \f 1/2 < \delta \le 1 $. Let \f 0 < c < d < \infty \A 
     and let \f J_{\pm}(c, d) \A be as in \r{\em (5.15)}$. Let \f u \A be 
     given by \r{\em (7.15)}$ with \f z \in J_{+}(c, d) \cup J_{-}(c, d) $. 
     Then there exists a positive constant \f C = C(c, d, \mu, \delta) \A 
     such that }
$$
          \M u \M_{1, -\delta} \le C \M f \M_{\delta}.           \T (7.23)
$$
\BP

        {\it Proof. \ } Suppose that (7.23) is not true. Then, for each
     \f n = 1, 2, \cdots $, there exist \f f_n \in L_{2, \delta}(\Rb) \A
     and \f z_n  \in J_{+}(c, d) \cup J_{-}(c, d) \A such that
$$
\left\{ \split
               \M f_n \M_{\delta} < 1/n, \\ *[4pt]
               \M u_n \M_{1, -\delta} = 1,

\endsplit \right.                                               \T (7.24)
$$
     where \f u_n = R(z_n)f_n $. Here we may assume with no loss of 
     generality that \f z_n \A converges to \f z_0 \A which is in the 
     closure of \f J_{+}(c, d) \cup J_{-}(c, d) $. We consider the case 
     that \f z_0 = \lambda \in [c, d] \A since the other case is much 
     easier. Then, using the Rellich selection theorem and the equation
$$
       -\, \Delta u_n - \mu(x)z_nu_n = \mu(x)f_n,                \T (7.25)
$$
     we see there exists a subsequence of \f \{ u_n \} \A which is a Cauchy
     sequence in \f H^1(\Rb)_{{\rm loc}} $. For the sake of simplicity of      
     notation we denote the subsequence by \f \{ u_n \} \A again. In view
     of Theorem 7.4, \f \{ u_n \} \A is a Cauchy sequence in 
     \f H_{-\delta}^1(\Rb) $. Let \f u_0 \A be the limit function.
     We have \f \M u_0 \M_{1, -\delta} = 1 $. It is easy to see that 
     \f u_0 \A is a (weak) solution of the homogeneous equation
     \f  -\, \Delta u_0 - \mu(x)\lambda u_0 = 0 $, and hence we have
     \f u_0 \in H_{-\delta}^2(\Rb) $. On the other hand, it follows from 
     Theorem 7.2 that, for \f n = 1, 2, \cdots $,
$$
\split
     & {\dm \M \CD u_n \M_{\delta-1, E_1}  
                          \le C\big(\M f_n \M_{\delta} 
                            + \M u_n \M_{-\delta} \big)  } \\ *[4pt]
     & \hspace{2.3cm} {\dm \le C( 1/n + 1 ) \le 2, }                          
\endsplit                                              \T (7.26)
$$
     and hence, by letting \f n \to \infty $, we have
     \f \M \CD u_0 \M_{\delta-1, E_1} < \infty $, i.e., \f u_0 \A satisfies
     the radiation condition. Therefore, by the two dimensional counterpart 
     of Corollary 3.8, we have \f u_0 = 0$, which contradicts the fact that
     \f \M u_0 \M_{1, -\delta} = 1 $. This completes the proof. \ \ $\Vert$            
\BP

        Let \f u = R(z)f \A be as in Theorem 7.4. Then it follows from 
     (7.23) and the equation \f -\, \Delta u - \mu(x)zu = \mu(x)f \A that
$$
          \M u \M_{2, -\delta} \le C \M f \M_{\delta}           \T (7.27)
$$
     with \f C = C(c, d, \mu, \delta) $. As in \S5, the operator norm in 
     \f \bB(H_s^j(\Rb), \, H_t^{\ell}(\Rb)) $ will be denoted by 
     \f \M \  \M_{(j,s)}^{(\ell,t)} $, where \f j, \ell = 0, 1, 
     2 $, \f s, t \in \Ro $. Thus we have 
\BP

        {\bf Theorem 7.6.} {\it Suppose that \r{\em Assumption 2.1}$ holds
     with \f N = 2 $. Let \f 1/2 < \delta \le 1 $. Let \f R(z) \A be the 
     resolvent of \f H $. Let \f 0 < c < d < \infty \A and let 
     \f J_{\pm}(c, \, d) \A be as above. Then there is a positive constant 
     \f C = C(c, d, \mu, \delta) \A such that
$$
         \M R(z) \M_{(0,\delta)}^{(2,-\delta)} \le C              \T (7.28)    
$$             
     for \f z \in J_+(c, \, d)\cup J_-(c, \, d) $.}       
\BP

        Now we can proceed as in \S6 to obtain the limiting absorption 
     principle for \f H \A with \f N = 2 $.
\BP

        {\bf Theorem 7.7.} {\it Suppose that \r{\em Assumption 2.1}$ holds
     with \f N = 2 $. Let \f 1/2 < \delta \le 1 \A. Let \f R(z) \A be the 
     resolvent of \f H $.
\SP

        \r{\em (i)}$ Let \f \lambda > 0 $. Then the extended resolvent
     \f R_{\pm}(\lambda) \A is well-defined by
$$
   R_{\pm}(\lambda) = \lim_{\eta\downarrow 0} R(\lambda \pm i\eta) \T (7.29)
$$
     in \f \bB(L_{2, \delta}(\Rb), \, H_{-\delta}^{2}(\Rb)) $.
\SP

       \r{\em (ii)}$ Let \f D_{\pm} \A be given by {\rm (6.6)} and extend 
     \f R(z) \A on \f D_{+} \A as in {\rm (6.7)}, i.e.,
$$
     R(\lambda+i\eta) = 
\left\{ \split
    R(\lambda+i\eta) \qq (\lambda > 0, \, \eta > 0), \\  *[4pt]      
    R_+(\lambda) \ \ \qq \ \ (\lambda > 0, \, \eta = 0).     
\endsplit \right.                                               \T (7.30)
$$
     Extend \f R(z) \A on \f D_{-} \A as in {\rm (6.8)}. Then
     \f R(z) \A is a 
     \f \bB(L_{2, \delta}(\Rb), \, H_{-\delta}^{2}(\Rb))$-valued
     continuoius function on each of \f D_+ \A and \f D_- $.
\SP

        \r{\em (iii)}$ For any \f z \in D_+ $ $\mb{\em [}$ or \f D_- $
     $\mb{\em ]}$, \f R(z) \A is a compact operator from 
     \linebreak \f L_{2,\delta}(\Rb) \A into \f H_{-\delta}^1(\Rb) $.      
\SP

        \r{\em (iv)}$ \  For \f 0 < c < d < \infty \A there exists a 
     constant \f C = C(c, d, \delta, m_0, M_0) > 0 \A such that, for 
     \f z \in \overline{J}_+(c, \, d) \cup \overline{J}_-(c, \, d) $,}
$$   
\left\{ \split
  {\dm \int_{E_s}(1 + r)^{-2\delta} 
                \big( |\gr R(z)f|^2 + |k|^2|R(z)f|^2 \big) \, dx } 
                                                              \\  *[5pt]
  \hspace{3.5cm} {\dm \le C^2 (1 + s)^{-(2\delta-1)} \M f \M_{\delta}^2 }
                                                              \\  *[4pt]
       \ \ \ \ \ \ \ \ \ \ \ \ \ \ \ \ \ \ \ \ \ \ \ \ \ \ \ \ \ \ \ \ \
         \ \ \ \ \ \ \ \ {\dm (s \ge 1, \ f \in L_{2,\delta}(\RN)), } 
                                                            \\  *[4pt] 
  {\dm \M \CD R(z)f \M_{\delta-1} \le C \M f \M_{\delta} \q
             ( f \in L_{2,\delta}(\RN)). }     
\endsplit \right.                                              \T (7.31)
$$
\SP
\BP

        {\bf Remark 7.8.} In the case that \f N \ge 3 $, all the constants
     \f C \A which appear in the evaluation of \f R(z) \A are constructive, 
     i.e., these constants \f C = C(c, d, \mu, \delta, \cdots ) \A can be 
     computed explicitly when the values of \f c $, \f d $, \f \mu $, 
     \f \delta $ \f \cdots \A are given. On the other hand, the constant 
     \f C \A in Proposition 7.5 is not constructive in our method, and
     hence the constant \f C \A in Theorem 7.6 is not constructive, too. 

$$
\ \ \ \ \ \ 
$$    
            
\cen{  {\bf References}  }           
\BP

\F   [1] S. Agmon, Spectral properties of Schr\"{o}dinger operators and 
     scattering theory, Ann. Scuola Sup. Pisa {\bf 2} (1975), 151-218. 
\SP

\F   [2] M. Ben-Artzi, Y. Dermanjian and J.-C. Guillot, Acoustic waves in
     perturbed stratified fluids: a spectral theory, Commun. Partial
     Differential Equations {\bf 14} (1989), 479-517.  
\SP

\F   [3] A. Boutet de Monvel-Berthier and D, Manda, Spectral and scattering 
     theory for wave propagation in perturbed stratified media, 
     Universit\"{a}t Bielefeld, BiBos, preprint, Nr. 606/11/93.
\SP

\F   [4] S. DeBi\'{e}vre and D. W. Pravica, Spectral analysis for optical 
     fibres and stratified fluids I: The liming absorption principle, J. 
     Functional Analysis {\bf 98} (1991) 406-436.
\SP

\F   [5] S. DeBi\'{e}vre and D. W. Pravica, Spectral analysis for optical 
     fibres and stratified fluids II: Absence of eigenvalues, Commun. Partial 
     Differential Equations {\bf 17}, (1992), 69-97.   
\SP

\F   [6] D. Eidus, The limiting absorption and amplitude problems for the 
     diffraction problem with two unbounded media, Comm. Math. Phys. 
     {\bf 107} (1986), 29-38.
\SP

\F   [7] T. Ikebe and Y. Sait\={o}, Limiting absorption method and absolute 
     continuity for the Schr\"odinger operator, J. Math. Kyoto Univ. {\bf 12} 
     (1972), 513-542.
\SP
             
\F   [8] W. J\"{a}ger and Y. Sait\={o}, The limiting absorption principle
     for the reduced wave operator with cylindrical discontinuity, Preprint.
     IWR (SFB359) 94-74, University of Heidelberg. 1994.
\SP

\F   [9] W. J\"{a}ger and Y. Sait\={o}, The limiting absorption principle
     for the reduced wave operator with multimedia. Preprint.1995.
\SP

\F   [10] G. Roach and B. Zhang, On Sommerfeld radiation conditions for the 
     diffraction problem with two unbounded media, Proc. Royal Soc. Edinburgh 
     {\bf 121A} (1992), 149-161.
\SP
 
\F   [11] Y. Sait\={o}, The principle of limiting absorption for second-order 
     differential equations with operator-valued coefficients, Publ. RIMS, 
     Kyoto Univ. {\bf 7} (1972), 518-619.
\SP

\F   [12] Y. Sait\={o}, The principle of limiting absorption for the 
     non-selfadjoint Schr\"o- \linebreak dinger operators in \f \RN \A 
     (\f N \neq 2 $), Publ. RIMS, Kyoto Univ. {\bf 9} (1974), 397-428.
\SP

\F   [13] Y. Sait\={o}, A remark on the limiting absorption principle for the 
     reduced wave equation with two unbounded media
     Pacific J. Math. {\bf 136} (1989), 183-208.

\SP

\F   [14] R. Weder, Absence of eigenvalues of the acoustic propagators in
     deformed waveguides, Rocky Mountain J. Math. {\bf 18} (1988), 495-503.
\SP

\F   [15] R. Weder, {\it Spectral and Scattering Theory for Wave Propagation 
     in Perturbed Stratified Media}, Springer-Verlaga, Berlin, 1991.
\SP

\F   [16] C. Wilcox, {\it Sound Propagation in Stratified Fluids}, 
     Springer-Verlag, New York, 1984.
\SP

\F   [17] B. Zhang, On radiation conditions for acoustic propagators in 
     perturbed stratified fluids. Preprint. 1994. 
          .
\end{document}